
\input amstex
\documentstyle{amsppt}
\magnification=1200
\hsize=16.5truecm
\vsize=23.3truecm

\catcode`\@=11
\redefine\logo@{}
\catcode`\@=13

\define\vth{\vartheta}

\define\ml{\operatorname{Exp-Lift}}

\define \bz{\Bbb Z}
\define \bq{\Bbb Q}
\define \br{\Bbb R}
\define \bc{\Bbb C}
\define \bh{\Bbb H}

\define\gi{\Gamma_{\infty}}
\define\gm{\Gamma}
\define\pd#1#2{\dfrac{\partial#1}{\partial#2}}

\define \Mn{S^nM}

\TagsOnRight
\document

\topmatter
\title
Elliptic genus of Calabi--Yau manifolds and 
Jacobi and Siegel modular forms
\endtitle
\author
V.  Gritsenko 
\footnote{Supported by Max-Planck-Institut
 f\"ur Mathematik in Bonn\hfill\hfill}
\endauthor

\address
St. Petersburg Department of Steklov Mathematical Institute,
\newline
${}\hskip 8pt $ Fontanka 27, 191011 St. Petersburg, Russia
\endaddress
\address
D\'epartement de Math\'ematique Universit\'e Lille I
\newline
${}\hskip 8pt $ 59655 Villeneuve d'Asq Cedex, France
\endaddress

\abstract
In the paper we study two types  of relations: a one
is between the elliptic genus of Calabi--Yau manifolds and 
Jacobi modular forms,
another one  is between the second quantized elliptic genus, Siegel 
modular forms and Lorentzian Kac--Moody Lie algebras.
We also  determine the structure of  the graded ring of 
the weak Jacobi forms with integral Fourier coefficients.
It gives us a number of applications to the theory of
elliptic genus and of the second quantized elliptic genus.
\endabstract
\rightheadtext
{Elliptic  genus of Calabi--Yau manifolds}
\leftheadtext{V. Gritsenko}

\endtopmatter

\document
\document

\head
Introduction 
\endhead

For  a   compact complex manifold one can define its elliptic genus
as a function in two complex variables.
If the first Chern class $c_1(M)$ of the complex manifold is equal to zero
in $H^2(M,\br)$, then the elliptic genus is an automorphic form in
variables $\tau\in \bh_1$ ($\bh_1$ is the upper-half plane) and $z\in \bc$.
More exactly, it is  a Jacobi modular form
with integral Fourier coefficients of weight $0$ and index $d/2$,
where $d=$dim${}_\bc(M)$. 
The $q^0$-term of the Fourier expansion ($q=e^{2\pi i \tau}$) 
of the elliptic genus is essentially equal to the Hirzebruch 
$\chi_y$-genus.
Thus we can analyze the arithmetic properties of the  
$\chi_y$-genus of the Calabi--Yau manifolds and 
its special values  such as signature ($y=1$) and Euler number ($y=-1$)
in  terms of    Jacobi modular forms.
The famous Rokhlin theorem about divisibility by $16$ of 
the signature  of a compact, oriented, differentiable spin manifold 
of dimension $4$ was one of the starting points of 
the theory of elliptic genera. 
Ochanine  generalized this Rokhlin
result to the manifolds of
$\hbox{dim}_\br M\equiv 4\mod 8$.
One can find an elegant proof the Ochanine's theorem 
using modular forms in one variable with respect to $\Gamma_0(2)$
in the lectures of Hirzebruch \cite{HBJ, Chapter 8}.

In this paper we  study the $\bz$-structure of the graded ring
$J_{0,*}^\bz=\oplus_{m\ge 1} J_{0,m}^\bz$ 
of all weak Jacobi forms with integral Fourier coefficients  
(Theorem 1.9).
We prove that this ring has four generators
$$
J_{0,*}^\bz=\bz[\phi_{0,1}, \phi_{0,2}, \phi_{0,3}, \phi_{0,4}]
$$
where $\phi_{0,1}, \dots, \phi_{0,4}$ are some fundamental Jacobi
forms related to Calabi--Yau manifolds of dimension $d=2,\,3,\,4,\,8$.
We consider some applications of this result to  Calabi--Yau manifolds.
Properties of the signature 
(i.e. the value of  $\chi_y$-genus  at $y=1$)  
modulo some powers of $2$ are well known (see \thetag{2.8}).
We analyze properties  of the value of the elliptic genus
and the  Hirzebruch $\chi_y$-genus at $y=-1$, 
$y=\zeta_3=e^{2\pi i/3}$ and $y=i$. 
For example, we prove that  the   Euler number
of a Calabi-Yau manifold  $M_d$ of dimension $d$ satisfies 
$$
e(M_d)\equiv 0\mod 8\qquad \text{if}
\quad d\equiv 2 \mod 8
$$
(see Proposition 2.4) and 
$$
\chi_{y=\zeta_3}(M_{d})\equiv 0\mod 9\qquad \text{if}
\quad d\equiv 2 \mod 6
$$
(see Proposition 2.5 for details).

The same Jacobi forms   $\phi_{0,1}, \dots, \phi_{0,4}$
are  generating functions for the multiplicities
of all positive roots of the four generalized Lorentzian Kac--Moody
Lie algebras of Borcherds type constructed in \cite{GN1--GN4}.
In \S 3 we show that the second quantized elliptic genus
(introduced by R. Dijkgraaf, E. Verlinde and
H. Verlinde  in \cite{DVV})
of an arbitrary Calabi--Yau manifolds of dimension $4$ and $6$ 
can be written as a product of some powers of  the denominator functions
of the  Lorentzian Kac--Moody algebras  constructed in \cite{GN4}.

The results of this paper were presented at the  seminar of Prof. Y. Manin 
and at the seminar "Arithmetic Lunch" of MPI in Bonn 
in the spring semester
of 1997. The author is  grateful to
Max-Planck-Institut f\"ur Mathematik  for hospitality.

\head
\S 1. Elliptic genus and weak Jacobi forms
\endhead

Let $M$ be a compact complex  manifold of complex dimension $d$
and $T_M$ be the holomorphic tangent bundle of  $M$.
One defines the formal power series 
$$
\Bbb E_{q,y} = \bigotimes_{n\geq 0}
{\bigwedge}_{-y^{-1}q^{n}}T_M
\otimes 
\bigotimes_{n\geq 1}{\bigwedge}_{-yq^n} {T}_M^*
\otimes 
\bigotimes_{n\geq 1}
S_{q^n} (T_M \oplus {T}_M^*)
$$
where $\wedge^k$ and $S^k$ denote the $k$-th exterior product and 
$k$-th symmetric 
tensor product  respectively and $T^*_M$ denotes the dual bundle. 
The elliptic genus of the complex manifold is the holomorphic 
Euler characteristic of this formal power series with  vector 
bundle coefficients.
This definition is standard in physical literature
(see, for example, \cite{W}, \cite{EOTY}, \cite{KYY}, \cite{AYS},
\cite{D}).
\definition{Definition}
The {\it elliptic genus} of $M$ is defined as follows 
$$
\chi(M;\,\tau,z):=y^{\frac{d}2}\int_M  \hbox{ch}(\Bbb E_{q,y})\, \hbox {td}(M)
\tag{1.1}
$$
where $\hbox{td}(M)$ is the Todd class of $M$. 
One applies  the Chern character
in $\hbox{ch}(\Bbb E_{q,y})$  
to all coefficients of $\Bbb E_{q,y}$
and the integral $\int_M$ denotes the evaluation of the top 
degree differential form on the fundamental cycle of the manifold.
\enddefinition
One can consider this  definition as a limit case of 
the level $N$ elliptic genus of Hirzebruch (see \cite{H}, \cite{HBJ},
\cite{K}, \cite{H\"o}).
See also \S 3 for a physical definition.
The coefficient $f(m,l)$ of the elliptic genus 
$$
\chi(M;\,\tau,z)=\sum_{n\ge 0, l} f(n,l) q^n y^l
$$
is equal
to the index of the Dirac operator
 (see, for example, \cite{HBJ, Appendix II}) twisted with the vector bundle 
$E_{m,l-\frac d{2}}$, where 
$
\Bbb E_{q,y} = \bigoplus_{m,\,n}  E_{m, n} 
\cdot q^m y^{n}
$.
In particular all coefficients of the elliptic genus  are integral.
According to  the Riemann--Roch--Hirzebruch theorem one can see
that the $q^0$-term of 
$\chi(M;\, \tau,z)$   is essentially
the Hirzebruch $\chi_y$-genus of the manifold $M$:
$$
\gather
\chi(M;\, \tau,z)=\sum_{p=0}^d (-1)^p\chi_p(M)\, y^{\frac{d}2-p}+
\tag{1.2}
\\
q\biggl(\sum_{p=-1}^{d+1}
(-1)^{p}y^{-p}
\bigl(\chi_p(M, T_M^*)-\chi_{p-1}(M, T_M^*)+
\chi_p(M, T_M)-\chi_{p+1}(M, T_M)\bigr)\biggr)
+\dots
\endgather
$$
where 
$
\chi(M, E)=\sum_{q=0}^d (-1)^q \,\hbox{dim}\, H^q (M, E)
$ 
and 
$
\chi^p(M, E)=\chi(M,\, \wedge^p T_M^*\otimes E)
$
or, for a K\"ahler manifold,
$\chi_p(M)=\sum_{q}(-1)^{q}h^{p,q}(M)$.

It turns out that the elliptic genus of a Calabi--Yau $d$-fold
is a weak Jacobi form of weight zero and index $d/2$. 
The Jacobi modular forms are modular forms
with respect to the Jacobi group or in other words with respect
to a parabolic subgroup of the integral
symplectic group.
In this paper we consider the Jacobi forms with respect to 
a maximal parabolic subgroup $\gi\subset Sp_2(\bz)$,
which consists of all element preserving an isotropic line.
The Jacobi group
$
\Gamma^J:=\gi/\{\pm E_4\} \cong SL_2(\bz)\ltimes H(\bz).
$
is the semi-direct product of $SL_2(\bz)$
and  the integral Heisenberg group $H(\bz)$
which is 
the central extension of $\bz\times \bz$.
The binary character
$
v_H([\lambda,\mu;\,\kappa]):=(-1)^{\lambda+\mu+\lambda\mu+\kappa}
$
can be extended to a binary character $v_J$ of the Jacobi group
if we put $v_J|_{SL_2(\bz)}=1$. 
In this paper we  use   Jacobi forms of integral 
or half integral index.

\definition{Definition}
A holomorphic  function
$\phi(\tau ,z)$ on $\bh_1\times \bc$ is called
a  {\it weak Jacobi form of
weight $k \in \bz/2$ and index $t \in \bz/2$ 
with respect to the Jacobi group $\Gamma^J$} if the function
$
\widetilde{\phi}(Z):=\phi(\tau ,z)\exp{(2\pi i t\omega )}
$
on the Siegel upper half-plane
$
\bh_2=\{Z=
\left(\smallmatrix \tau&z\\z&\omega\endsmallmatrix\right)
\in M_4(\bc), \ \hbox{Im} (Z)>0\}
$
is 
a $\Gamma^J$-modular form of weight $k$
with  character (or multiplier system) $v_J^{2t}$,
i.e., if it satisfies the functional equations
$$
\phi(\frac{a\tau+b}{c\tau+d},\,\frac z{c\tau+d})=
(c\tau+d)^k e^{2\pi i t c z^2/(c\tau+d)}\phi(\tau,z)
\qquad (\left(\smallmatrix a&b\\c&d\endsmallmatrix\right) \in SL_2(\bz))
\tag{1.3a}
$$
and 
$$
\phi(\tau, z+\lambda \tau+ \mu)=
(-1)^{2t(\lambda+\mu)} 
e^{-2\pi i t(\lambda^2 \tau+2\lambda z)}\phi(\tau,z)
\qquad (\lambda, \mu\in \bz).
\tag{1.3b}
$$
and it has the Fourier expansion of the type
$$
\phi(\tau ,z)=\sum\Sb n\ge 0,\ l \in t+\bz\\
\vspace{0.5\jot} \endSb
f(n,l) \exp{(2\pi i(n\tau +lz))}.
$$
We denote the space of all weak Jacobi forms of weight $k$
and index $t$ by $J_{k,t}$.
\enddefinition

If the Fourier coefficients satisfy  $f(n,l)=0$ unless $4tn-l^2\ge 0$,
then $\phi$ is holomorphic  at infinity.
We call such Jacobi forms  {\it holomorphic} Jacobi forms.
Weak Jacobi forms of integral index were firstly defined 
in \cite{EZ}. The main advantage of this notion is that
$J_{k,t}$ is still finite dimensional and  
the  graded ring of all weak Jacobi forms is finitely generated
(over $\bc$).
In \cite{EZ} only Jacobi forms of integral index were considered.
In the next example we define the main week Jacobi forms 
of half-integral index (see Lemma 1.4 bellow).
\example{Example 1.1} 
{\it Weak Jacobi forms $\phi_{0, \frac{3}2}$ and 
$\phi_{-1, \frac{1}2}$}.
Let us define a  Jacobi theta-series
$$
\align
\vartheta(\tau ,z)&=\hskip-2pt\sum\Sb n\equiv 1\, mod\, 2 \endSb
\,(-1)^{\frac{n-1}2}
\exp{(\frac{\pi i\, n^2}{4} \tau +\pi i\,n z)}=
\sum_{m\in \bz}\,\biggl(\frac{-4}{m}\biggr)\, q^{{m^2}/8}\,y^{{m}/2}\\
{}&=
-q^{1/8}y^{-1/2}\prod_{n\ge 1}\,(1-q^{n-1} y)(1-q^n y^{-1})(1-q^n)
\tag{1.4}
\endalign
$$
where we use the formal variables $q=e^{2\pi i \tau}$, $y=e^{2\pi i z}$.
This is a holomorphic Jacobi form of weight $1/2$ and index $1/2$ with 
multiplier system $v_\eta^3\times v_H$
($v_\eta$ is the multipier system of the Dedekind eta-function 
$\eta(\tau)$).
One can  check that 
$$
\align
{}&\phi_{0, \frac 3{2}}(\tau,z)=
\hskip -2pt \frac{\vartheta(\tau ,2z)}{\vartheta(\tau ,z)}
=
y^{-\frac 1{2}}
\prod_{n\ge 1}(1+q^{n-1}y)(1+q^{n}y^{-1})(1-q^{2n-1}y^2)
(1-q^{2n-1}y^{-2})\hskip -1pt\in \hskip -2pt J_{0,\frac 3{2}}\\
{}&\phi_{-1, \frac 1{2}}(\tau,z)=\frac{\vartheta(\tau ,z)}{\eta(\tau)^3}
=-y^{-1/2}\prod_{n\ge 1}\,(1-q^{n-1} y)(1-q^n y^{-1})(1-q^n)^{-2}
\in J_{-1,\frac 1{2}}
\endalign
$$
are weak Jacobi forms with integral Fourier coefficients.
\endexample

The next proposition is  well-known in physical literature
(see \cite{KMY}, \cite{AYS}). See also  \cite{H\"o}
where another normalization
for the elliptic genus function was used.

\proclaim{Proposition 1.2}
If $M_d$ is a compact complex  manifold of dimension $d$ with $c_1(M)=0$
(over $\br$), 
then its elliptic genus 
$\chi(M_d;\, q,y)$ is a weak Jacobi form of weight $0$ and index 
$d/2$.
\endproclaim
\demo{Proof} For the convenience of the readers we give here 
a proof  which is similar to the proof of 
the modular behavior of the level $N$ elliptic genus (see 
\cite{H} and \cite{HBJ, Appendix III}).
Let us represent $\chi(M;\tau, z)$ 
in terms of the theta-series. 
Let $c_i(T_M)\in H^{2i}(M, \bz)$ be the Chern class  of $T_M$
and let  $c(T_M)$ and $\hbox{ch}(T_M)$ be the total Chern class
and the Chern character of $M$ and
$x_i=2\pi i \xi_i$ ($1\le i\le d$)
be the formal
Chern parameters of $T_M$:
$$
c(T_M)=\sum_{i=0}^r c_i(T_M)=\prod_{i=1}^{r}(1+x_i), \qquad
\hbox{ch}(T_M)=\sum_{i=1}^r e^{x_i}.
$$
We recall that
$$
\hbox{ch\,}({\bigwedge}_t E)=\prod_{i=1}^r(1+te^{x_i}),
\qquad
\hbox{ch\,}(S_t E)=\prod_{i=1}^r \frac{1}{1-te^{x_i}}.
$$
Then we have
$$
\hbox{ch\,}({\Bbb E}_{q,y})\,\hbox{td\,}(T_M)=
\prod_{i=1}^d\, \prod_{n=1}^{\infty}
\frac{(1-q^{n-1}y^{-1}e^{-x_i})(1-q^{n}ye^{x_i})}
{(1-q^{n-1}e^{-x_i})(1-q^{n}ye^{x_i})}x_i\,.
$$
Therefore
$$
\Phi(\tau, z, \xi_1, \dots, \xi_d)=y^{d/2} 
\hbox{ch\,}({\bold E}_{q,y})\hbox{td\,}(T_M)=
\prod_{i=1}^d \frac{\vth(\tau, z+\xi_i)}
{\vth(\tau, \xi_i)} (2\pi i \xi_i)
$$
is holomorphic function in $\tau$, $z$ and $\xi_i$.
(We recall that $\vth(\tau, z)$ has zero of order $1$ along $z=0$.)
For arbitrary 
$\gamma=\left(\smallmatrix a&b\\c&d\endsmallmatrix\right) \in SL_2(\bz)$
one has 
$$
\frac{\vth(\gamma<\tau>, \frac{z+\xi_i}{c\tau+d})}
{\vth(\gamma<\tau>,\  \frac{\xi_i}{c\tau+d})}=
e^{\pi i  c (z^2+ 2z \xi_i)/(c\tau+d)}
\frac{\vth(\tau,\  {z+\xi_i})}{\vth(\tau, \xi_i)}.
$$
Thus $\Phi(\tau, z, \xi_1, \dots, \xi_d)$ transforms like 
a Jacobi form of weight $d$ (due to $d$ factors $\xi_i$ after the quotient
of the theta-series) and 
index $\frac{d}2$ if  
$c_1(M)=x_1+\dots x_d=0$ (over $\br$!).
Therefore the coefficient of arbitrary  monom
$\xi_1^{n_1}\cdot \dots \cdot \xi_d^{n_d}$ of total degree $d$ 
in the Taylor expansion
of $\Phi(\tau, z, \xi_1, \dots, \xi_d)$ with respect to 
$\xi_i$ is a  weak  Jacobi form of weight $0$ and index $\frac{d}2$.
It follows that after integration over $M$ one gets 
a Jacobi form of weight zero and index $d/2$.
The claim that the elliptic genus does not contain 
Fourier coefficients with negative powers of $q$ 
follows directly from the definition \thetag{1.1}.
\enddemo

\remark{Remark} The definition \thetag{1.1} has sense and Proposition 1.2
is still true for arbitrary compact manifold with a $Spin^c$-structure. 
Moreover one can give a similar definition
for an arbitrary complex vector bundle $E$ with $c_1(E)=0$.
A generalization of this staff for  vector bundles with 
$c_1(E)\ne 0$ see in \cite{G2}.
\endremark 

\example{Example 1.3} {\it Elliptic genus of Enriques and  $K3$ surfaces.} 
Let us consider
the case of a complex surface $M_2$ with $c_1(M_2)=0$. 
Then
$$
\chi(M_2;\, \tau, z)\in J_{0,1}=\bc\, \phi_{0,1}(\tau,z)
$$
where 
$$
\phi_{0,1}(\tau,z)
=(y+10+y^{-1})+q(10y^{-2}-64y^{-1}+108-64y+10y^2)+q^2(\dots )
\tag{1.5}
$$
is the unique, up to a constant, weak Jacobi form of weight $0$
and index $1$.
The value of arbitrary weak Jacobi form of weight $0$ at $z=0$ 
is a holomorphic
$SL_2(\bz)$-modular form of weight $0$, i.e. it is a constant.
According \thetag{1.2}
$$
\chi(M_d;\tau,0)=\sum_{p=0}^d (-1)^p\chi_p(M)=e(M_d)
\tag{1.6}
$$
is the Euler number of $M_d$. Therefore we have 
$$
\chi(\hbox{Enriques} ;\, \tau, z)=\phi_{0,1}(\tau,z), \quad
\chi(K3;\, \tau, z)=2\phi_{0,1}(\tau,z)
$$
for arbitrary Enriques and  $K3$ surfaces.
\endexample
The structure of the bigraded algebra 
over $\bc$ of all Jacobi forms of integral index was determined 
in \cite{EZ}. We describe the structure of
the graded $\bz$-algebra of all Jacobi forms of  weight zero
with integral Fourier coefficients.
Firstly we show how the weak  Jacobi forms of integral 
and half-integral indices  are related. 
\proclaim{Lemma 1.4} Let  $m$ be integral, then we have for the 
weak Jacobi forms
$$
J_{2k,m+\frac 1{2}}=\phi_{0, \frac 3{2}}\cdot J_{2k,m-1},
\qquad
J_{2k+1,m+\frac 1{2}}=\phi_{-1, \frac 1{2}}\cdot J_{2k+2,m}
$$
where $\phi_{0, \frac 3{2}}$ and $\phi_{-1, \frac 1{2}}$
are defined in Example 1.1.
\endproclaim
\demo{Proof}
From the functional equation \thetag{1.3b} one has 
$$
\hbox{div}(\phi_{2k,m+\frac 1{2}})
\supset \{z\equiv \frac 1{2}, \frac \tau{2}, \frac {\tau+1}{2} 
\mod  \bz\tau+\bz\},
\qquad
\hbox{div}(\phi_{2k+1,m+\frac 1{2}})
\supset \{\bz\tau+\bz\}.
$$
Thus any Jacobi form of half-integral index and even (resp. odd) weight
is divisible by $\phi_{0, \frac {3}2}$ (resp.  $\phi_{-1, \frac {1}2}$).
\enddemo

\example{Example 1.5} {\it Calabi--Yau $3$-folds and $5$-folds}.
Let $c_1(M_3)=0$ and $c_1(M_5)=0$.
From Lemma 1.4 it follows that 
$$
J_{0,\frac 3{2}}=\bc\, \phi_{0, \frac {3}2},\qquad
J_{0,\frac 4{2}}=\bc\, \phi_{0, \frac {3}3}\phi_{0,1}.
$$
Thus 
$$
\align
\chi(M_3;\, q,y)&
=\frac{e(M_3)}2\,
\bigl[(y^{1/2}+y^{-1/2})+q(-y^{5/2}+y^{1/2}+y^{-1/2}-y^{-5/2})+q^2(\dots)
\bigr],\\
\chi(M_5;\, q,y)&
=\frac{e(M_5)}{24}
\bigl[(y^{\pm \frac 3{2}}+11y^{\pm \frac 1{2}})+
q(-y^{\pm \frac 7{2}}+\dots)+q^2(-11y^{\pm \frac9{2}}+\dots)\bigr]
\endalign
$$
($y^{\pm l}$ means that we have two summands with
$y^{l}$ and $y^{-l}$ respectively).
As a simple corollary  we have that
{\it the Euler characteristic of an arbitrary Calabi--Yau
$5$-fold is divisible by  $\,24$} and 
$$
\chi_1(CY_5)=\frac{1}{24}\, e(CY_5), \qquad
\chi_2(CY_5)=\frac{11}{24}\, e(CY_5).
$$
In particular, for the  Hodge numbers of an arbitrary  strict Calabi--Yau 
$5$-fold we have 
$$
11(h^{1,1}-h^{1,4})=h^{2,2}-h^{2,3}+10(h^{2,1}-h^{3,1}).
$$
\endexample
To calculate the elliptic genus of a  Calabi--Yau manifold
 of even dimension
we introduce some other basic Jacobi forms.
Let us denote by 
$J^\bz_{0,m}$ the $\bz$-module of all  weak Jacobi forms
of weight $0$ and index $m$ with integral Fourier  coefficients.
We consider the graded ring
$$
J^\bz_{0,*}=\bigoplus_{m\in \bz_{\ge 0}}J^\bz_{0,m}
$$
of all Jacobi forms of integral index with integral Fourier
coefficients 
and its ideal
$$
J^\bz_{0,*}(q)=\{\phi \in J^\bz_{0,*}\,|\, 
\phi(\tau,z)=\sum_{n\ge 1, \ l\in \bz} a(n,l) q^n y^l\}
$$
of the  Jacobi forms without $q^0$-term.

\proclaim{Lemma 1.6} The ideal $J^\bz_{0,*}(q)$ is  principal.
It is  generated by a weak Jacobi form of weight $0$  and index $6$
$$
\xi_{0,6}(\tau,z)=\Delta(\tau)\phi_{-1,\frac {1}2}(\tau,z)^{12}=
\frac{\vartheta(\tau,z)^{12}}{\eta(\tau)^{12}}=
q(y^{\frac{1}2}-y^{-\frac{1}2})^{12}+q^2(\dots).
$$
\endproclaim
\demo{Proof} Let $\psi\in J_{k,m}$ be an arbitrary Jacobi form.
The product 
$$
\Psi(\tau, z)=\exp\bigl(-8\pi^2m\, G_2(\tau)z^2 \bigr)
\psi(\tau, z)
$$
transforms like a $SL_2(\bz)$-modular form of weight $k$.
We recall that 
$
G_2(\tau)=-\frac 1{24}+\sum_{n=1}^{\infty}\sigma_1(n)\,q^n
$ 
is a quasi-modular form of weight $2$, i.e.
it satisfies
$$
G_2(\gamma<\tau>)=(c\tau+d)^2G_2(\tau)-\frac {c(c\tau+d)}{4\pi i}.
$$
Thus   we have
$
\Psi(\frac{a\tau+b}{c\tau+d},\,\frac z{c\tau+d})=
(c\tau+d)^k
\Psi(\tau,z)
$.
Therefore the coefficients $f_n(\tau)$
in the Taylor expansion 
$$
\Psi(\tau,z)=\sum_{n\in \bz} f_n(\tau) z^n
\tag{1.7}
$$
are  $SL_2(\bz)$-modular forms of weight $k+n$.
If $\psi\in J_{0,m}^{\bz}(q)$ has no $q^0$-term, then 
$f_{2n}(\tau)$ is a cusp form of weight $2n$. 
Thus $n\ge 6$ and 
$\psi(\tau,z)$ has zero of order at least $12$ along $z=0$.
Therefore $\psi(\tau,z)/\xi_{0,6}(\tau,z)$ is holomorphic
on $\bh_1\times \bc$.
\enddemo

\proclaim{Corollary 1.7} A weak Jacobi form of weight $0$ is uniquely
determined by its $q^0$-term if its index is less than $6$ or equal
to $\frac {13}{2}$.
\endproclaim

\example{Example 1.8} {\it Calabi--Yau  manifolds of dimension 
$4$, $6$, $8$  and $10$.}
The last corollary shows us that for these values of dimension 
of a Calabi--Yau manifold $M_d$
one can define the  elliptic genus of $M_d$ if one knows 
its Hirzebruch $\chi_y$-genus. 
In \cite{GN1--GN4} we used the following important
Jacobi forms of small index
$$
\align
\phi_{0,2}(\tau ,z)&=
{\tsize\frac{1}2} \eta(\tau )^{-4}
\sum_{m\,,n\in \bz}
{(3m-n)}\biggl(\frac{-4}m\biggl)\biggl(\frac{12}n\biggl)
q^{\frac{3m^2+n^2}{24}}y^{\frac{m+n}2}\\
{}&=(y+4+y^{-1})+q(y^{\pm 3}-8y^{\pm 2}-y^{\pm 1}+16)
+q^2(\dots),
\tag{1.8}\\
\vspace{1\jot}
\phi_{0,3}(\tau ,z)&=\phi_{0,\frac{3}2}^2(\tau ,z)=
(y+2+y^{-1})+q(-2y^{\pm 3}-2y^{\pm 2}+2y^{\pm 1}+4)+q^2(\dots),\\
\vspace{1\jot}
\phi_{0,4}(\tau ,z)&=
\frac {\vartheta(\tau ,3z)}{\vartheta(\tau ,z)}
=(y+1+y^{-1})-q(y^{\pm 4}+y^{\pm 3}-y^{\pm 1}-2)
+q^2(\dots).
\tag{1.9}
\endalign
$$
Using these Jacobi forms together with  $\phi_{0,1}$
(see \thetag{1.5}) one can construct a  basic of the module 
of the 
weak  Jacobi forms of weight $0$  and index
$2$, $3$, $4$ and $5$ and to write a formula for the elliptic genus
of $M_d$ in terms of the cohomological invariants $\chi_p(M_d)$:
$$
\align
\psi_{0,2}^{(1)}&=\phi_{0,2}=y+4+y^{-1}+q(\dots)\\
\psi_{0,2}^{(2)}&=\phi_{0,1}^2-20\phi_{0,2}=(y^2+22+y^{-2})+q(\dots)\\
\vspace{2\jot}
\psi_{0,3}^{(1)}&=\phi_{0,3}=\phi_{0,\frac 3{2}}^2=y+2+y^{-1}+q(\dots)\\
\psi_{0,3}^{(2)}&=\phi_{0,1}\phi_{0,2}-15\phi_{0,3}
=(y^{\pm 2}-y^{\pm 1}+12)+q(\dots)
\tag{1.10}\\
\psi_{0,3}^{(3)}&=
\phi_{0,1}^3-30\phi_{0,1}\phi_{0,2}+117\psi_{0,3}=
\psi_{0,1}|T_-(3)-3\phi_{0,1}
=(y^{\pm 3}+34)+q(\dots)\\
\vspace{2\jot}
\psi_{0,4}^{(1)}&=\phi_{0,4}=y+1+y^{-1}+q(\dots)\\
\psi_{0,4}^{(2)}&=\phi_{0,1}\phi_{0,3}-12\phi_{0,4}
=(y^{\pm 2}+10)+q(\dots) \tag {1.11} \\
\psi_{0,4}^{(3)}&=\phi_{0,2}\psi_{0,2}^{(2)}-4\psi_{0,4}^{(2)}-24\phi_{0,4}
=(y^{\pm 3}-y^{\pm 1}+24)+q(\dots)\\
\psi_{0,4}^{(4)}&=\psi_{0,2}^{(2)}|T_-(2)-2\phi_{0,4}^{(2)}=
(y^{\pm 4}+46)+q(\dots),
\endalign
$$
where $T_-(m)$ is  the standard  Hecke operator
$\bigl(\phi_{0,t}|T_-(m)\bigr)(\tau,z)=\sum_{n,l}f_m(n,l)q^ny^l$
with 
$f_m(N,L)=m\sum_{a|(N,L,m)} a^{-1}f(\frac{Nm}{a^2},\frac{L}a)$.

Fourier coefficient $f(n,l)$ of a weak Jacobi form of weight $0$ 
and index $t$ depends only on the norm $4nt-l^2$ of its index
and $\pm l \mod 2t$.
For every Jacobi form written above its $q^0$-term contains
all ``orbits" of non-zero
Fourier coefficients with negative norm $4nt-l^2$ of its index.
This fact will be important in \S 3 bellow.
For $d=10$ the situation is a little bit  more complicated. 
After some calculation we obtain  the following  Jacobi forms
$$
\align
\psi_{0,5}^{(1)}&=5y+2+5y^{-1}+q(-y^{5}+\dots)\\
\psi_{0,5}^{(2)}&=(y^{\pm 2}+y^{\pm 1}+8)+q(y^5+\dots)  \\
\psi_{0,5}^{(3)}&=(y^{\pm 3}+y^{\pm 1}+20)+q(\dots)
\tag {1.12}
\\
\psi_{0,5}^{(4)}&=(y^{\pm 4}-y^{\pm 1}+36)+q(-3y^{5}+\dots)\\
\psi_{0,5}^{(5)}&=(y^{\pm 5}+58)+q(\dots)
\endalign
$$
where we  include in the formulae all ``orbits" of non-zero
Fourier coefficients with  negative norm $20n-l^2$.
One knows (see \cite{EZ}) that 
$\hbox{dim\,} J_{0,m}=m$ for 
$m=1$, $2$, $3$, $4$ and $5$. Thus 
the functions $\psi_{0,m}^{(p)}$ constructed above
($1\le p\le m$) form a basis of 
$J_{0,m}^{\bz}$. Thus for $d=2,\ 4,\ 6,\ 8,\ 10$ we can represent 
the elliptic genus of a compact complex manifold $M_d$ with 
$c_1(M_d)=0$ as a linear combination
$$
\chi(M_d;\tau,z)=\sum_{p=0}^{\frac d{2}-2} 
\chi_p(M_d)\psi_{0,\frac d{2}}^{(\frac d{2}-p)}(\tau,z)+
c \psi_{0,\frac d{2}}^{(1)}(\tau,z)
$$
with $c\in \bz$.
In particular,  the coefficients $\chi_p$ of the Hirzebruch $y$-genus
of the manifold $M_d$ 
satisfy
the following relations
$$
\aligned
\chi_2&=22\chi_0-4\chi_1\\
\chi_3&=34\chi_0 -14\chi_1+2\chi_2\\ 
\chi_4&=46\chi_0-25\chi_1+10\chi_2-\chi_3
\endaligned
\qquad
\aligned
{}&\text{and }\qquad e(M_4)\equiv 0 \mod 6\ (\text{ for CY}_4)\\
{}&\text{and }\qquad e(M_6)\equiv 0 \mod 4\ (\text{ for CY}_6)\\ 
{}&\text{and }\qquad e(M_8)\equiv 0 \mod 3\ (\text{ for CY}_8).
\endaligned
\tag {1.13}
$$
For a Calabi--Yau $10$-fold we have
$$
\chi_5=58\chi_0-36\chi_1+20\chi_2-8\chi_3+
\frac{2}5 (\chi_4+\chi_3-\chi_2-\chi_1).
$$
We note that for dimension $d\ge 12$ ($d\ne 13$) an additional  
term $c\xi_{0,6}$
could appear in the formula for the elliptic genus.
The odd dimensions $d=7$, $9$, $11$ and $13$
are related to the even  dimensions considered above
by the factor $\phi_{0, 3/2}$. For example,
$$
\chi(M_7; \tau,z)=\chi_2(M_7)\phi_{0, \frac 3{2}}\phi_{0,2}
-\chi_1(M_7)\phi_{0, \frac 3{2}}\psi_{0, 2}^{(2)},
\qquad 
e(M_7)=12(\chi_2-4\chi_1)\equiv 0\mod 12.
$$
\endexample

These examples  show us that information about 
Jacobi forms provides some information about Calabi--Yau manifolds.
In the next theorem we construct a basis of   the module $J^\bz_{0,m}$
and we find generators of the graded ring $J_{0,*}$.
\proclaim{Theorem 1.9} {\bf 1}. Let $m$ be  a positive integer.
The module
$$
J^\bz_{0,m}/ J^\bz_{0,m}(q)=
\bz\,[\psi_{0,m}^{(1)},\dots ,\psi_{0,m}^{(m)}]
$$
is a free $\bz$-module of rank $m$. Moreover there is a basis
consists of 
$\psi_{0,m}^{(n)}$ ($1\le n\le m$)
with the following $q^0$-term
$$
\align
[\psi_{0,m}^{(n)}]_{q^0}&=y^n+\dots+ y^{-n}
\qquad\qquad\qquad (3\le n\le m)\\
[\psi_{0,m}^{(2)}]_{q^0}&=y^2-4y+6-4y^{-1}+y^{-2}\\
[\psi_{0,m}^{(1)}]_{q^0}&=\frac{1}{(12,m)}\,
\bigl(my+(12-2m)+my^{-1}\bigr)
\endalign
$$
where $(12,m)$ is the greatest common divisor of $12$ and $m$.

{\bf 2}. The graded ring of all  weak Jacobi forms
of weight $0$ with integral coefficients is finitely generated
$$
J^\bz_{0,*}=
\bigoplus_{m}J^\bz_{0,m}=
\bz\,[\phi_{0,1}, \phi_{0, 2}, \phi_{0,3}, \phi_{0,4}]
$$
where the corresponding  Jacobi forms are defined in 
\thetag{1.8}--\thetag{1.9}.
$\phi_{0,1}$, $\phi_{0, 2}$, $\phi_{0,3}$
are algebraicly independent and 
$$
4\phi_{0,4}=\phi_{0,1}\phi_{0, {3}}-\phi_{0,2}^2.
$$
\endproclaim
To prove Theorem 1.9 we need
\proclaim{Lemma 1.10}For arbitrary weak Jacobi
form  $\phi_{0,m}=\sum_{n,l}f(n,l)q^ny^l$ 
the following  identities are valid
$$
m\sum_{l}f(0,l)=6\sum_{l}l^2 f(0,l), \qquad 
24m\sum_{l} f(0,l)=\sum_{n}(m-6n^2)f(1,n).
$$
\endproclaim
\demo{Proof}
We consider the Taylor expansion \thetag{1.7}
$$
\exp\bigl(-8\pi^2m\, G_2(\tau)z^2 \bigr)
\phi_{0,m}(\tau, z)=
f_0+\sum_{n\ge 1} f_n(\tau) z^n
$$
where  $f_n(\tau)$ is a  $SL_2(\bz)$-modular form
of weight $n$. Thus $f_2(\tau)\equiv 0$.
The left hand side of the identities of the lemma are 
equal to  the first two coefficients of $f_2(\tau)$.
This proves the lemma.
We remark that one can continue the list of similar identities
calculating the second, third, $\dots$ coefficients of 
$f_2(\tau)$. 
\enddemo
Using \thetag{1.2} we get 
$$
\frac {e(M)d}{12}=\sum_{p}(-1)^p\chi^p(M)(\frac{d}2-p)^2
\tag{1.14}
$$ 
if $c_1(M)=0$ over $\br$.
\remark{Remark}Identity \thetag{1.14} was rediscovered several
times in the mathematical and physical literature.
See, for example, \cite{LW}, \cite{S} for a proof based on
the Riemann--Roch--Hirzebruch formula and \cite{AYS} where 
\thetag{1.14} is related to the sum rule for the charges
in $N=2$ super-conformal field theory.
Lemma 1.11 provides us with an automorphic proof of \thetag{1.14}
(see \cite{GN3, Lemma 2.2} for a more general statement).
\endremark
\demo{Proof of Theorem 1.9}
For a weak Jacobi form  the polynom $[\phi]_{q^0}$  is of order not greater
than $m$ (see \cite{EZ}). It can not be
a constant according to the first identity of Lemma 1.10.
Moreover if   $[\phi]_{q^0}$ is ``linear", i.e.
$[\phi]_{q^0}=ay+b+ay^{-1}$, then $a$ is divisible by 
$\frac{m}{(m,12)}$. It follows that $m$ Jacobi forms of type 
dicribed in  the first statement  of the theorem, if they would exist,
form a basis of the module $J^\bz_{0,m}/ J^\bz_{0,m}(q)$.

Let us construct a basis by induction.
For $m=1$, $2$, $3$, $4$, $6$, $8$, $12$ we set $\psi_{0,m}^{(1)}=\phi_{0,m}$
where $\phi_{0,1}$, $\dots$, $\phi_{0,4}$ were defined above  
and 
$$
\align
\phi_{0,6}(\tau,z)&=\phi_{0,2}\phi_{0,4}-\phi_{0,3}^2=
(y+y^{-1})+q(\dots),\\
\phi_{0,8}(\tau,z)&=\phi_{0,2}\phi_{0,6}-\phi_{0,4}^2
=(2y-1+2y^{-1})+q(\dots),\\
\phi_{0,12}(\tau,z)&=\phi_{0,4}\phi_{0,8}-2\phi_{0,6}^2
=(y-1+y^{-1})+q(\dots).
\endalign
$$
Now we can define the weak Jacobi forms $\psi_{0,m}^{(1)}$ 
using the following procedure. (We omit the weight $0$ bellow
to simplify the notation). Let  $m\ge 5$. We set
$$
\psi_{m,I}=\widetilde{\psi}_{m-4}\phi_4+\widetilde{\psi}_{m-2}\phi_2
-2\widetilde{\psi}_{m-3}\phi_3=my^{\pm 1}+(12-2m)+q(\dots)
$$
where by definition  $\widetilde{\psi}_{m}:=(12,m){\psi}_{m}^{(1)}$.
We can take $\psi_{m,I}$ as    ${\psi}_{m}^{(1)}$ 
if  $m$ is coprime with $12$.
If $m\equiv 0\mod 2$ and $m\ge 6$ we can define a form 
$\psi_{m,II}=\frac{1}{2}\psi_{m,I}$ with integral Fourier
coefficients.        
If $m\equiv 0\mod 3$ and $m\ge 9$ we define
$$
\psi_{m,III}=\frac{2}{3}\widetilde{\psi}_{m-3}\phi_3
+\frac{1}{3}\widetilde{\psi}_{m-6}\phi_6
-\widetilde{\psi}_{m-4}\phi_4
=\frac{m}{3}y^{\pm 1}+(4-\frac{2m}3)+q(\dots).
$$
For $m\equiv 0\mod 4$ and $m> 12$ we take 
$$
\psi_{m,IV}=\frac{1}{4}[\widetilde{\psi}_{m-12}\phi_{12}
+\widetilde{\psi}_{m-4}\phi_4
-\widetilde{\psi}_{m-8}\phi_8]
=\frac{m}{4}y^{\pm 1}+(3-\frac{m}2)+q(\dots).
$$
In the remain  cases $m\equiv 0\mod 6$ and $m\equiv 0\mod 12$
($m>12$) we put
$$
\align
\psi_{m,VI}&=\frac{1}{2}\psi_{m,III}
=\frac{m}{6}y^{\pm 1}+(2-\frac{m}3)+q(\dots),\\
\psi_{m,XII}&=\frac{2}{3}\widetilde{\psi}_{m-3}\phi_3
-\frac{1}2\widetilde{\psi}_{m-4}\phi_4
-\frac{1}{6}\widetilde{\psi}_{m-6}\phi_6
+\frac{1}{12}\widetilde{\psi}_{m-12}\phi_{12}
=\frac{m}{12}y^{\pm 1}+\frac{6-m}6+q(\dots).
\endalign
$$
Thus we   finish the construction of functions ${\psi}_{m}^{(1)}$ if 
we put 
$
{\psi}_{m}^{(1)}:=\psi_{m,D}
$
for  $(m,12)=D$.

Next we construct ${\psi}_{m}^{(2)}$. We put
$$
\psi_{0,m}^{(2)}=\widetilde{\psi}_{0,m-3}\phi_{0,3}
-\widetilde{\psi}_{0,m-4}\phi_{0,4}
-\widetilde{\psi}_{0,m}\qquad (m\ge 5)
$$
and
$$
\psi_{0,2}^{(2)}=\phi_{0,1}^2-24\phi_{0,2},\ \ 
\psi_{0,3}^{(2)}=\phi_{0,1}\phi_{0,2}-18\phi_{0,3},\ \ 
\psi_{0,4}^{(2)}=\phi_{0,1}\phi_{0,3}-16\phi_{0,4}.
$$
One can see that these forms have the $q^0$-term equal
to $y^2-4y+6-4y^{-1}+y^{-2}$.
For $3\le n\le m-2$ one can use $\phi_{0,3}$ and
the forms $\phi_{0,m-3}^{(n-1)}$ in order to construct
$\phi_{0,m}^{(n)}$. To finish the proof of the first statement of 
the theorem we put 
$$
\phi_{0,m}^{(m-1)}=\phi_{0,1}^{m-2}\cdot \phi_{0,2},\qquad
\phi_{0,m}^{(m)}=\phi_{0,1}^{m}.
$$
We note that 
each function of the basis $\{\psi_{m}^{(n)}\}_{n=1}^m$ of the module
$J^\bz_{0,m}/ J^\bz_{0,m}(q)$ constructed above is 
a polynom in the basic  Jacobi forms
$\phi_{0,1}$, $\dots$, $\phi_{0,4}$.
The torsion relation 
$4\phi_{0,4}=\phi_{0,1}\phi_{0,3}-\phi_{0,2}^2$ follows from Corollary 1.7.
The form $\xi_{0,6}$ generating the principle ideal 
$J^\bz_{0,*}(q)$ is also a polynom in 
$\phi_{0,1},\dots, \phi_{0,4}$:
$$
\xi_{0,6}=-\phi_{0,1}^2\phi_{0,4}
+9\phi_{0,1}\phi_{0,2}\phi_{0,3}-8\phi_{0,2}^3-27\phi_{0,3}^2.
\tag{1.15}
$$
(To prove the last formula one needs to  check  that the $q$-constant term
of the write hand side is zero and to compair coefficients
at the first power $q$.)

Let us prove that $\phi_{0,1}$, $\phi_{0,2}$ and $\phi_{0,3}$ are algebraicly
independent. For this end we consider its values at $z=\frac{1}{2}$.
We have
$$
\phi_{0,2}(\tau, \frac 1{2})\equiv 2,\quad
\phi_{0,3}(\tau, \frac 1{2})\equiv 0, \quad
\phi_{0,4}(\tau, \frac 1{2})\equiv -1.
$$
(The two  last identities follow from definition and the first
one is a corollary of the torsion relation.) 
The restriction of 
$$
\phi_{0,1}(\tau, \frac 1{2})
=\alpha(\tau)=8+2^8q+2^{11}q^2+11\cdot 2^{10}q^3+
3\cdot 2^{14}q^4+359\cdot 2^9 q^5+\dots
\tag{1.16}
$$
is a modular function with respect to $\Gamma_0(2)$ with  a character
of order $2$ (see \S 2).
The square of this  function is, up to factor $2^{12}$, the ``Hauptmodul"
for the congruence subgroup group $\Gamma_0(2)$.
Since only one function obtained from the  Jacobi forms of different indices
is a non-constant function for $z=1/2$, then 
they are algebraicly independent.
\enddemo

\head
\S 2. Special values of the elliptic genus
\endhead

In this section we analyze the value  of the elliptic genus 
at the following special points
$z=0$ (Euler number), $z=\frac 1{2}$ (signature),
$z=\frac{\tau+1}2$ ($\hat A$-genus) and $z=\frac 1{3},
\ \frac 1{4},\ \frac 1{6}$.
For this end we have to study the restriction of the main
generators of the graded ring of the integral week Jacobi forms.
A special value of a Jacobi form
is a modular form  in $\tau$.
In the next lemma we give a little more
precise statement than in  \cite{EZ, Theorem 1.3}.
\proclaim{Lemma 2.1} Let $\phi\in J_{0,t}$ ($t\in \bz/2$) and 
$X=(\lambda, \mu)\in \bq^2$. Then
$$
\phi|_{X}(\tau,0)=\phi(\tau, \lambda\tau+\mu)
\exp(2\pi i t(\lambda^2\tau+ \lambda \mu))
$$ 
is an automorphic form of weight $0$ with a character 
with respect to  the  subgroup 
$$
\Gamma_X=\{M\in SL_2(\bz)\,|\, XM-X\in \bz^2\}.
$$ 
\endproclaim
Let us consider some particular examples which we need in this section.

\example{Example 2.2} Let $X=(0, \frac{1}{N})$.
Then 
$$
\Gamma_X=\Gamma_0^{(1)}(N)=
\biggl\{M=\pmatrix a&b\\c&d\endpmatrix
 \in SL_2(\bz)\,|\, M\equiv  \pmatrix 1&*\\0&1\endpmatrix \mod N
\biggr\}.
$$
If $t$ is integer, then 
the character of the function 
$$
\phi|_X(\tau,0)=\phi(\tau,\frac{1}N)
$$
is given by $v(M)=\exp(-2\pi i tc/ N^2)$ and it has
has order $\frac N{(t,N)}$.
\endexample

It is easy to see that if $\phi\in J_{k,m}^\bz$ with integral $m$, then 
the form $\phi(\tau,\frac 1{N})$ still has integral Fourier coefficients if
$N=1, \dots, 6$. In particular, the value of $\xi_6(\tau,z)$ at 
these points is related to the ``Hauptmodule" for the corresponding
group $\Gamma_0(N)$:
$$
\aligned
\xi_6(\tau,\frac{1}2)&=2^{12}\frac{\Delta(2\tau)}{\Delta(\tau)},\\
\xi_6(\tau,\frac{1}3)&=
3^{6}\left(\frac{\Delta(3\tau)}{\Delta(\tau)}\right)^{1/2},
\endaligned
\qquad
\aligned
\xi_6(\tau,\frac{1}4)&=
2^{6}\left(\frac{\Delta(4\tau)}{\Delta(\tau)}\right)^{1/2},
\\
\xi_6(\tau,\frac{1}6)&=
\left(\frac{\Delta(\tau)\Delta(6\tau)}
{\Delta(2\tau)\Delta(3\tau)}\right)^{1/2}.
\endaligned
\tag{2.1}
$$
Let us analyze the corresponding values of the four generators $\phi_{0,n}$
of the graded ring $J_{0,*}^{\bz}$. From the definition 
(see \thetag{1.8}--\thetag{1.9}) and the identity
$4\phi_{0,4}=\phi_{0,1}\phi_{0, {3}}-\phi_{0,2}^2$
we obtain 
$$
\phi_{0,1}(\tau, 0)=12,\quad 
\phi_{0,2}(\tau, 0)=6,\quad
\phi_{0,3}(\tau, 0)=4,\quad
\phi_{0,4}(\tau, 0)=3
\tag{2.2}
$$
and
$$
\aligned
\phi_{0,1}(\tau, \frac{1}2)&=\alpha(\tau)\\
\phi_{0,2}(\tau, \frac{1}2)&=2\\
\phi_{0,3}(\tau, \frac{1}2)&=0\\
\phi_{0,4}(\tau, \frac{1}2)&=-1
\endaligned
\qquad
\aligned
\phi_{0,1}(\tau, \frac{1}3)&=\beta^2(\tau)\\
\phi_{0,2}(\tau, \frac{1}3)&=\beta(\tau)\\
\phi_{0,3}(\tau, \frac{1}3)&=1\\
\phi_{0,4}(\tau, \frac{1}3)&=0
\endaligned
\qquad
\aligned
\phi_{0,1}(\tau, \frac{1}4)&=\frac{\gamma(\tau)^4+4}{\gamma(\tau)}\\
\phi_{0,2}(\tau, \frac{1}4)&=4\gamma^2(\tau)\\
\phi_{0,3}(\tau, \frac{1}4)&=2\gamma(\tau)\\
\phi_{0,4}(\tau, \frac{1}4)&=1.
\endaligned
\tag{2.3}
$$
The automorphic functions $\alpha(\tau)$, $\beta(\tau)$ and $\gamma(\tau)$
are automorphic forms of weight $0$ with respect to the group
$\Gamma_0(2)$, $\Gamma_0^{(1)}(3)$ and $\Gamma_0^{(1)}(4)$ 
respectively.
These functions have integral Fourier coefficients.
The identity \thetag{1.15} gives us the following relations
for the modular forms $\alpha$, $\beta$ and $\gamma$
$$
\gather
2^{12}\frac{\Delta(2\tau)}{\Delta(\tau)}=
\alpha(\tau)^2-64, \qquad
3^{6}\left(\frac{\Delta(3\tau)}{\Delta(\tau)}\right)^{1/2}
=\beta(\tau)^3-27\\
2^{6}\left(\frac{\Delta(4\tau)}{\Delta(\tau)}\right)^{1/2}
=4(\left(\frac{\gamma(\tau)}2\right)^2
-\left(\frac 2{\gamma(\tau)}\right)^2).
\endgather
$$
It follows that
$$
\alpha(\tau)-8\equiv 0 \mod 2^8,\qquad
\beta(\tau)-3\equiv 0 \mod 3^3
\tag{2.4}
$$
(compare with \thetag{1.16}).
Using the  definition of $\phi_{0,3}$ and $\gamma(\tau)$ 
and the relations between 
the Jacobi theta-series $\vartheta_{ab}$ of level $2$  we have 
$$
\gamma(\tau)=
\frac{\vartheta_{00}(2\tau)}{\vartheta_{01}(2\tau)}=
\frac{\vartheta_{00}(2\tau,0)}{\vartheta_{01}(2\tau,0)}.
\tag{2.5}
$$
Using  Corollary 1.7 we check that
$$
\phi_{0,1}(\tau, 2z)=\phi_{0,2}^2(\tau,z)-8\phi_{0,4}(\tau,z).
$$
Thus
$$
\alpha(\tau)=16\gamma(\tau)^4-8=
16\frac{\vth_{00}^4(2\tau)}{\vth_{01}^4(2\tau)}-8.
\tag{2.6}
$$
In particular {\it all Fourier coefficients of $\gamma(\tau)$
and $\alpha(\tau)$ are positive}.

In connection with \thetag{2.5}--\thetag{2.6} we note that 
one can write the generators  of the graded ring  
as  symmetric polynomials  
in  
$\xi_{ab}(\tau,z)=\vartheta_{ab}(\tau,z)/\vartheta_{ab}(\tau,0)$:
$$
\gather
\phi_{0,1}=4(\xi_{00}^2+\xi_{10}^2+\xi_{01}^2),
\qquad
\phi_{0,\frac 3{2}}=
4 \xi_{00}\xi_{10}\xi_{01}
\\
\vspace{2\jot}
\phi_{0,2}=2\bigl((\xi_{00}\xi_{10})^2+
(\xi_{00}\xi_{01})^2+(\xi_{10}\xi_{01})^2\bigl).
\tag{2.7}
\endgather
$$
To check these formulae one can use Corollary 1.7 and the fact that
the generators of modular group transform $\xi_{ab}$ to each other.

\example{Example 2.3}
Let $X=(\frac{1}{N}, \frac{1}{N})$.
Then  $\Gamma_X$ contains the principle congruence subgroup
$\Gamma_1(N)$. In some cases $\Gamma_X$ will be strictly  larger.
For example, if $X_2=(\frac{1}{2}, \frac{1}{2})$,  then
$$
\phi|_{X_2}(\tau,0)=\phi(\tau,\frac{\tau+1}2)
\exp(\frac{\pi i}2({\tau+1}))
$$
is an automorphic form with respect of the so-called theta-group
$$
\Gamma_\theta=
\biggl\{ M=\pmatrix a&b\\c&d\endpmatrix \in 
 SL_2(\bz)\,|\  M\equiv \pmatrix 1&0\\0&1\endpmatrix 
\ \text{ or }\ \pmatrix 0&1\\1&0\endpmatrix
\mod 2
\biggr\}.
$$
The corresponding character is given by
$\epsilon_2(M)=\exp(2\pi i m(d+b-a-c)/4)=\pm 1$. 
This  character is trivial if the index $m$ of Jacobi
form is even and it has order $2$ if the index is odd.
Let us  consider $\Gamma_\theta$-automorphic  function
$$
\hat \phi_{m}(\tau)=
q^{-\frac{m}4}\phi_{0,m}(\tau, -\frac{\tau+1}2).
$$
We have
$$
\hat\phi_3= 0,\quad \hat\phi_4= -1,\quad
\hat\phi_2= -2,\quad
\hat\xi_6=\hat \phi_1^2+64=\left(\frac{\vth_{00}}{\eta}\right)^{12}
$$
where  
$$
\hat\phi_1(\tau)=4\frac{\vth_{10}^4-\vth_{01}^4}{\vth_{01}^2\vth_{10}^2}
=-q^{-\frac 1{4}}+20q^{\frac 1{4}}+\dots 
\in \frak M_0^\bz(\Gamma_\theta, \epsilon_2).
$$
\endexample

Let us analyze some special values  of the elliptic genus.
As it easy follows from \thetag{1.2} we get Euler number and
signature for $z=0$ ($d$ is arbitrary) and $z=\frac{1}2$ 
($d$ is even)
$$
\align
\chi(M_d, \tau, 0)&=e(M_d),\\
\chi(M_d, \tau, \frac{1}2)&=\sigma_M(\tau)=
(-1)^{\frac{d}2}s(M_d)+q(\dots)\in 
\frak M_0^\bz(\Gamma_{0}(2), v_2),
\quad
v_2(\left(\smallmatrix a&b\\c&d\endsmallmatrix\right))=
e^{\pi i m \frac {c}2}.
\endalign
$$
In \S 1 we calculated  Euler number $e(M_d)$ for small $d$
in terms of Jacobi forms. It gave us some divisibility
of Euler number of Calai--Yau manifolds.
We note that the quotient $e(M)/24$ appears in physics as obstruction
to cancelling the tadpole (see \cite{SVW} where it was proved
that $e(M_4)\equiv 0 \mod 6$).
\proclaim{Proposition 2.4}Let $M_d$ be an almost complex manifold of complex 
dimension $d$ such that $c_1(M)=0$ in $H^2(M, \br)$.
Then 
$$
d\cdot e(M_d)\equiv 0 \mod 24.
$$
If $c_1(M)=0$ in $H^2(M, \bz)$,  then we have a more strong congruence
$$
e(M)\equiv 0 \mod 8 \qquad\text{if }\  d\equiv 2 \mod 8.
$$
\endproclaim
\demo{Proof}The first fact follows simply from 
\thetag{2.2} or from \thetag{1.14}. If $d\equiv 2 \mod 8$
one can write the elliptic genus as a polynom over $\bz$ in 
the generators $\phi_{0,*}$
$$
e(M_d)\equiv P(\phi_{0,1}, \phi_{0,2}, \phi_{0,3}, \phi_{0,4})|_{z=0}
\equiv c_{1,m}(\phi_{0,1}|_{z=0})(\phi_{0,4}|_{z=0})^{\frac {d-2}8}
\mod 8.
$$
If one put $z=-\frac{\tau+1}2$, i.e.,
$y=-q^{1/2}$ (see Example 2.3), then one obtain that the series 
$$
\Bbb E_{q,-q^{1/2}} = \bigotimes_{n\geq 1}
{\bigwedge}_{q^{n/2}}T_M
\otimes 
\bigotimes_{n\geq 1}{\bigwedge}_{q^{n/2}} {T}_M^*
\otimes 
\bigotimes_{n\geq 1}
S_{q^n} (T_M \oplus {T}_M^*)
$$
is $*$-symmetric.  According to the  Serre duality all Fourier 
coefficients of 
$\hat\chi(M_d,\tau)$ are even.
The constant $ c_{1,m}$ from the last congruence
is equal to the coefficient of $\hat\chi(M_d,\tau)$ at 
the minimal negative power of $q$. Therefore
$c_{1,m}$ is even and we obtain divisibility of $e(M_{8m+2})$ by 8.
\enddemo
We note that divisibility of $d\cdot e(M)$ by $3$  was proved by 
F. Hirzebruch in 1960.
For  a hyper-K\"ahler compact  manifold the claim of the proposition above
was proved  by S. Salamon  in \cite{S}.
After my talk on the  elliptic genus  at a seminar of MPI
in Bonn in April 1997 Professor F. Hirzebruch informed me that 
the result of Proposition 2.4 was known for him (non-published).
Using some natural examples he also proved that 
this property of divisibility 
of the Euler number modulo $24$ is strict (see \cite{H2}).

Formulae \thetag{2.3} provide us with a formula for the signature 
$\chi(M_d;\tau, \frac{1}2)$ as a polynom in $\alpha(\tau)$.
As a corollary of \thetag{2.3} and Theorem 1.9
we have that 
for an arbitrary Jacobi form of integral index
$$
\aligned
\phi_{0,4m}(\tau,\frac{1}2)&=c+2^{13}q(\dots)\\
\phi_{0,4m+2}(\tau,\frac{1}2)&=2 c+2^{12}q(\dots)
\endaligned
\qquad 
\aligned 
\phi_{0,4m+1}(\tau,\frac{1}2)&=8 c+2^{8}q(\dots)\\
\phi_{0,4m+3}(\tau,\frac{1}2)&=16 c+2^{9}q(\dots).
\endaligned
\tag{2.8}
$$
Similar to the proof of Proposition 2.4 we obtain a better congruence
for the signature of a manifold with dim$\equiv 2 \mod 8$
and $c_1(M)=0$:
$$
\chi(M_{8m+2};\tau,z)=16c+2^{9}q(\dots). \tag{2.9}
$$
{\it This gives us  another proof of the Oshanine's theorem
in this particular case.}

This is interesting that the values of the Hirzebruch $y$-genus
at $y=e^{2\pi i/3}$ and $y=i$ also have some properties of divisibility.
For $z=\frac{1}3$ (resp. $z=\frac{1}4$) we can write 
$\phi_{0,m}(\tau, \frac{1}3)$ 
(resp. $\phi_{0,m}(\tau, \frac{1}4)$)
as a polynom in $\beta(\tau)=3+27(q+\dots)$ (resp. in $\gamma(\tau)^{\pm 1}$).
This gives us the following results
$$
\gather 
\phi_{0,3m}(\tau,\frac{1}3)=c+3^{6}q(\dots),\qquad
\phi_{0,3m+1}(\tau,\frac{1}3)=9c+3^{4}q(\dots)\\
\phi_{0,3m+2}(\tau,\frac{1}3)=3c+3^{3}q(\dots).
\tag{2.10}
\endgather
$$
Thus  we have 

\proclaim{Proposition 2.5} If $c_1(M)=0$ (over $\br$), then
$$
\gather
\chi(M_{6m}; \tau,\frac{1}3)=c_1 \mod 3^6,\qquad
\chi(M_{6m+2}; \tau,\frac{1}3)=9c_2 \mod 3^4,\\
\chi(M_{6m+4}; \tau,\frac{1}3)=3c_3 \mod 3^3.
\endgather
$$
\endproclaim

We finish with some relations for $z=\frac 1{4}$:
$$
\gather
\chi(M_{8m+2}; \tau,\frac{1}4)=4c+2^{4}q(\dots),\qquad
\phi_{0,4m+2}(\tau,\frac{1}4)=4c+2^{5}q(\dots)\\ 
\phi_{0,4m+3}(\tau,\frac{1}4)=2c+2^{8}q(\dots).
\tag{2.11}
\endgather
$$

\head
\S 3. Second quantized elliptic genus (SQEG) of Calabi--Yau manifolds  
\endhead

\subhead
3.1. SQEG of Calabi--Yau  manifolds and Siegel modular forms
\endsubhead
The notion of elliptic genus of $N=2$ super-symmetric theories 
was introduced more than $10$ years ago
(see, for example, \cite{W1}, \cite{W2},  \cite{EOTY}, 
\cite{AYS}, \cite{KYY}).
In physics the elliptic genus of a Calabi--Yau manifold
$M_d$
is defined as 
 the genus one partition function of the super-symmetric sigma 
model whose target space is   $M_d$.
By definition this is 
the trace of an operator 
over the Ramond-Ramond sector of 
the sigma model
$$
\chi(M_d;\, \tau, z)=\hbox{Tr}_{\Cal H}\, (-1)^F q^{L_0-\frac c{24}}
\overline{q}^{\overline{L}_0-\frac c{24}}y^{F_L}
$$
where 
$q^{L_0-\frac c{24}}$ is  the evolution operator,
$F_L$ is the fermion number and $L_0$ ($\overline{L_0}$)
is Virasoro generator of of left (right) movers, 
$(-1)^F$ is the fermion parity operator in $N=2$ super-symmetric theory
and $c$ is  the central charge. 

We can consider   $n$-fold symmetric product of the manifold  $M$,
i.e., the orbifold space $\Mn=M^n/S_n$, 
where $S_n$ is the symmetric group of $n$ elements.
This is a singular manifold but one can define the   {\it orbifold}
elliptic genus of $\Mn$ (see for details the talk 
of R. Dijkgraaf at ICM-1998 in Berlin \cite{D}).
It gives us a string version of the  elliptic genus of the symmetric 
product of a Calabi--Yau manifold.
One can compare this construction with the definition of the
string Euler number of symmetric products 
(see \cite{HH}).
Using some  arguments from the conformal field theory 
on orbifolds it was proved in \cite{DVV} and  \cite{DMVV} 
that the string elliptic genus
of the second quantization 
$\cup_{n\ge1} \Mn$ of  a Calabi--Yau manifold $M$
coincides with the second quantized elliptic genus of the 
given manifold:
$$
\sum_{n=0}^\infty p^n \chi_{orb}(\Mn;\, q,y)  
= \prod_{m\ge 0,\, l,\, n>0}
\frac 1{
(1- q^m y^l p^n)^{f(mn, l)}}
\tag{3.1}
$$
where 
$$
\chi(M,\tau,z)=\sum_{m\ge 0,\ l\in \bz \ (or\,\bz/2)} f(m,l)q^my^l
$$
is the elliptic genus of $M$.

For a $K3$ surface, the product in the left hand side of \thetag{3.1}
is essentially  the power $-2$  of the infinite product expansion 
of
the product of all even theta-constants
(see \cite{GN1})
$$
\Delta_5(Z)^2=2^{-12}\prod_{{}^tab=0\, mod\,2} \Theta_{a,b}(Z)^2=
(qpr) \prod_{(m,l,n)>0} (1- q^m y^l p^n)^{f(mn, l)}
$$
where 
$$
Z=\pmatrix \tau&z\\z&\omega\endpmatrix\in \bh_2,\quad
q=e^{2\pi i \tau},\ y=e^{2\pi i r},\ p=e^{2\pi i \omega}
$$
belongs to the Siegel upper half-plane,
$\Theta_{a,b}(Z)$ are the even Siegel theta-constants of level $2$
and $f(m,l)$ is the Fourier coefficient of 2$\phi_{0,1}(\tau,z)$
(see Example 1.3).
The modular form $\Delta_5(Z)$ is the first cusp form with respect 
to the Siegel modular group  $Sp_4(\bz)$ 
(with a non-trivial binary character).
In \cite{DVV} \thetag{3.1} was conjectured as
a generalization of a formula
for the  orbifold Euler numbers 
$$
\sum_{n\ge 1}\, e_{orb} (\Mn)\,p^n=\prod_{n>0}(1-p^n)^{-24}=
\prod_{m\ge 0,\, l, \,n>0}
(1- q^m y^l p^n)^{-f(mn, l)} |_{y=1}.
$$
In the case of an algebraic  $K3$ surface 
$e_{orb}(\Mn)$ is the topological Euler characteristic of 
the Hilbert scheme $M^{(n)}$ of zero dimension subschemas of length $n$
(see \cite{G\"o}).

\smallskip
Following \cite{DVV, \S 4}  we call the product
$$
\Psi(M;\, Z)=\prod_{m\ge 0,\, l,\, n>0}
(1- q^m y^l p^n)^{-f(mn, l)},
\tag{3.2}
$$
where $f(m,l)$ are coefficients of the elliptic genus of $M$,
the {\it second-quantized elliptic genus} (SQEG) of the manifold  $M$.

In this section we prove that the SQEG multiplied by 
a factor ({\it Hodge anomaly}) depending only 
on invariants $\chi_p(M)$
of the Calabi--Yau $d$-fold 
$$
E(CY_d;\, Z)= \text{Hodge anomaly}(CY_d;\, Z) \cdot \Psi (CY_d;\,Z)
$$
is a Siegel automorphic form
with respect to the paramodular group of polarization  $(1,d)$
(resp. $(1,2d)$) for even (resp. odd) dimension $d$.
We  also calculate $E(CY_d;\, Z)$ for small $d$ in terms 
of the basic Siegel modular forms for the corresponding paramodular groups.
For small dimensions $d=2$, $3$, $4$, $6$  the function
$E(CY_d;\, Z)$
is product of some powers of the denominator function
of generalized Lorentzian Kac--Moody Lie algebras constructed
in \cite{GN1--GN4}.

We considered in \cite{GN1}--\cite{GN4}    the following 
infinite product for an arbitrary nearly
holomorphic Jacobi form  
$\phi_{0,t}(\tau,z)=\sum_{n,l} c(n,l)q^ny^l$ of weight $0$ and integral
index $t$
$$
\ml(\phi_{0,t})(Z)=
q^{A}y^Bs^{C}
\prod\Sb n,l,m\in \bz\\
\vspace{0.5\jot}
(n,l,m)>0\endSb
 (1-q^ny^ls^{tm})^{c(nm,l)},
\tag{3.3}
$$
where
$$
A=\frac{1}{24}\sum_{l}f(0,l),\quad
B=\frac{1}{2}\sum_{l>0}lf(0,l),\quad
C=\frac{1}{4}\sum_{l}l^2f(0,l)
$$
and $(n,l,m)>0$ means that 
if $m> 0$, then $l$ and $n$ are arbitrary integers,
if $m=0$, then $n>0$ and $l\in \bz$
or  $l<0$ if $n=m=0$.

In order to define the corresponding product  for 
Calabi--Yau manifolds 
of odd dimension it is useful to modify
the definition of SQEG.
We remark that for  arbitrary $\phi_{0,\frac{d}2}(\tau, z)
\in J_{0,\frac d{2}}$
the Jacobi form  $\phi_{0,\frac d{2}}(\tau, 2z)\in J_{0,2d}$
has even index.
To see this one can apply  symplectic transformation 
$\Lambda_2=\hbox{\rm diag}(1,2,1,2^{-1})$ to
the corresponding  $\Gamma^J$-modular form:
$$
(F|\Lambda_2) 
\left(\smallmatrix \tau&z\\z&\omega\endsmallmatrix\right)=
F\left(\smallmatrix \tau& 2z\\ 2z& 4\omega\endsmallmatrix
\right).
$$
The next theorem is a particular case of Theorem 2.1 
from  \cite{GN4}.
\proclaim{Theorem 3.1} Let $M=M_d$ be a compact complex manifold
of dimension $d$
with trivial  $c_1(M)$,
$$
\chi(M;\, \tau,z)=\sum_{m\ge 0,\ l\in \bz \ (or\,\bz/2)} f(m,l)q^my^l
$$
be its elliptic genus and $\ \hbox{\rm SQEG}(M;\,Z)$  ($Z\in \bh_2$)
be its second quantized elliptic genus \thetag{3.3}.
We define a factor
$$
H(M;\,Z)=
\cases
\eta(\tau)^{-\frac 1{2}(e-3\chi'_{d_0})}
\prod_{p=1}^{d_0}\bigl(\vartheta(\tau, pz)\,e^{\pi i p^2\omega}\bigl)
^{-\chi'_{d_0-p}}&\ \text{ if }\  d=2d_0\\
\eta(\tau)^{-\frac 1{2}{e}}
\prod_{p=1}^{d_0}\bigl(\vartheta(\tau, \frac{2p-1}2 z)\,
e^{\frac {1 }4 \pi i (2p-1)^2\omega}\bigl)
^{-\chi'_{d_0-p+1}}&\ \text{ if } \ d=2d_0+1
\endcases
$$
where $e=e(M)$ is Euler number of $M$ and 
$\chi'_p=(-1)^p\chi_p(M)$  (see \thetag{1.2}).
Then  the product
$$\align
E(M;\,Z)&=\Psi(M;\,Z)\cdot \hbox{\rm SQEG}(M;\,Z)\qquad (d=2d_0)\\
E^{(2)}(M;\,Z)&=(E |\Lambda_2) (M;\,Z)\qquad \qquad (d=2d_0+1)
\endalign
$$
determines a Siegel automorphic form of weight $-\frac 1{2}{\chi'_{d_0}(M)}$
if $d$ is even and of of weight $0$ if $d$ is odd 
with a character  or a multiplier system 
of order ${24}/(24, e)$ 
with respect to a double extension of the  paramodular group 
$\Gamma_d^+$ (resp. $\Gamma_{2d}^+$), if $d$ is even
(resp. $d$ is odd).
The  divisor of $E(M;\,Z)$  (resp. of $E^{(2)}(M;\,Z)$)
on $\Cal A_d^+=\Gamma_d\setminus \bh_2$ (resp. $\Cal A_{2d}^+$)
is the union of  a finite number of Humbert modular surfaces 
$H_D(b)$  of discriminant $D=b^2-4ad$ (resp. $D=4b^2-8ad$)
with  multiplicities
$$
m_{D,b}=-\sum_{n>0}f(n^2a,nb).
$$
\endproclaim

\remark{Remark}We call the factor $H(M;\,Z)$ defined above  Hodge anomaly
of SQEG. The divisor of $\vartheta(\tau, pz)$ is $pz\in \tau\bz+\bz$.
Thus the zeros and poles of $H(M;\,Z)$ (but not of $E(M;\,Z)$!)
 are completely defined by  $\chi_p(M)$.
\endremark
We would like to consider applications of this result to
Calabi--Yau manifolds of dimension $3$, $5$ and $2$, $4$, \dots , $10$.
Before  doing this we recall the definitions of the notions we used above.
The paramodular group $\Gamma_t$ is isomorphic to the integral symplectic
group of the skew-symmetric form with elementary divisors
$(1,t)$. It can be realized as a subgroup
of $Sp_4(\bq)$
$$
\Gamma_t:=
\left\{\left(\smallmatrix
*  & t* & * & *\\
*  & * & * & t^{-1}*\\
*  & t * & * & *\\
t* & t*& t* & *
\endsmallmatrix\right)
 \in Sp_4 (\bq)\,|\, \text{ all $*$ are integral}\, \right\}.
$$
If $t\ne 1$, the group $\Gamma_t$ is not a maximal discrete subgroup
and it has normal extensions (see, for example, \cite{GH1}).
By definition,
$$
\Gamma_t^+=\Gamma_t\cup \Gamma_tV_t, \qquad
V_t={\tsize\frac 1{\sqrt{t}}}\left(\smallmatrix
0   &   t   &   0   &   0\\
1  & 0 &   0   &   0\\
0   &   0   &   0   &   1\\
0   &   0   &   t   &   0
\endsmallmatrix\right).
$$
The double quotient
$$
\Cal A_t\overset{2:1}\to\longrightarrow
\Cal A_t^+=\Gamma_t^+\setminus \bh_2
$$
of $\Cal A_t$ can be interpreted as a moduli space of
lattice-polarized $K3$ surfaces for arbitrary $t$ or as the moduli
space of Kummer surfaces of $(1,p)$-polarized Abelian surfaces for a
prime $t=p$ (see \cite{GH1, Theorem 1.5}).
Any Humbert surface  in $\Cal A_t^+$
of discriminant $D$
can be represented in the form
$$
H_D^+(b)=\pi_t^+\,\biggl(\bigcup
\Sb g\in \gm_t^+  \endSb
g^*(\{Z\in \bh_2\,|\,a\tau+bz+t\omega=0\})\biggr)
$$
where  $a,b\in \bz$, $D=b^2-4ta$,
$0\le b<2t$  and 
$\pi_t^+: \bh_2\to \Cal A_t^+$ is the natural projection.

According  \thetag{3.3} the infinite product $E(M;\,Z)$ can be written
as follows
$$
E(M_d;\,Z)=\cases 
\ml(-\chi(M_d;\,\tau,z))&\ \text{ if $d$ is even}\\
\ml(-\chi(M;\,\tau, 2z))&\ \text{ if $d$ is odd}.
\endcases
$$
We note that for odd $d$ the product 
$E(M;\,Z)$ (without $\Lambda_2$-modification) is an automorphic function
with respect to the group $\Lambda_2 \Gamma_{2d} \Lambda_2 ^{-1}$
conjugated to the paramodular group.

\subhead
3.2. SQEG of Calabi--Yau manifolds of even dimension 
and Lorentzian  Kac--Moody algebras
\endsubhead
\smallskip
{\bf 1. The case of $\bold C\bold Y_{\bold 2}$.}
One of the  starting points for  considerations
in \cite{DVV}
was the infinite product expansion formula for the modular form
$\Delta_5(Z)$.
According to Example 1.3 we have for the SQEG of surfaces
with trivial $c_1$
$$
\align
E(Enriques\  surface;\,Z)&=\ml(-\phi_{0,1})(Z)=\Delta_5(Z)^{-1},\\
E(K3;\,Z)&=\ml(-2\phi_{0,1})(Z)=\Delta_5(Z)^{-2}.
\endalign
$$
In \cite{GN1} we proved that the modular form $\Delta_5(Z)$ 
defines an automorphic correction of the Kac--Moody  
algebra of hyperbolic type with the Cartan matrix 
$$
A_{1,II}=\pmatrix \hphantom{-}2&-2&-2\\
-2&\hphantom{-}2&-2\\
-2&-2&\hphantom{-}2\endpmatrix.
$$
It means that the Fourier expansion of the  modular form  $\Delta_5(Z)$
determines a generalized Kac--Moody Lie super-algebra with a system 
of the real simple roots of type $A_{1,II}$.
We would like to note that the  Lorentzian Kac--Moody algebras 
with Cartan matrix of types $A_{1,0}$ and $A_{1,I}$  constructed
in \cite{GN2} are related to SQEG
of {\it some vector bundles of rank $2$ over a manifold of dimension}
$14$ (see \cite{G2}).

\smallskip
{\bf 2. The case of $\bold C\bold Y_{\bold 4}$.}
The basic  Jacobi modular forms for this dimension are 
the Jacobi forms $\phi_{0,2}$ and $\psi_{0,2}^{(2)}$
(see \thetag{1.8} and \thetag{1.10}). They correspond to the following
cusp forms for the paramodular group  $\Gamma_2$ 
(see \cite{GN1} and \cite{GN4}):
$$
\align
\Delta_2(Z)&=\ml(\phi_{0,2}(\tau,z))=
\hbox{Lift}(\eta(\tau)^3\vartheta(\tau,z))\\
{}&=\sum_{N\ge 1}\
\sum\Sb
 n,\,m >0,\,l\in \bz\\
\vspace{0.5\jot} n,\,m\equiv 1\,mod\,4\\
\vspace{0.5\jot} 2nm-l^2=N^2
\endSb
\hskip-4pt
N\biggl(\frac {-4}{Nl}\biggr)
\sum_{a\,|\,(n,l,m)} \biggl(\frac {-4}{a}\biggr)
\, q^{n/4} y^{l/2} s^{m/2}
\in \frak M_2^{cusp}(\Gamma_2,v_\eta^6\times v_H)
\endalign
$$
and
$$
\Delta_{11}(Z)=\hbox{Lift}(\eta(\tau )^{21}\vartheta(\tau ,2z))=
\ml(\psi_{0,2}^{(2)}(\tau,z))
\in \frak N_{11}(\Gamma_{2}).
$$
For an arbitrary  Calabi--Yau $4$-fold $M_4$
we have the following formula for its
SQEG
$$
E(M_4;\,Z)=\Delta_{11}(Z)^{-\chi_0(M)} \Delta_{2}(Z)^{\chi_1(M)}.
\tag{3.4}
$$
Its divisor is  equal to
$
(\chi_1(M)-\chi_0(M))H_1-\chi_{0}H_4
$.
We note that $\Delta_{2}(Z)^4$ is the first $\Gamma_2$-cusp form
with trivial character and $\Delta_{11}(Z)$ is the first cusp form of odd
weight with respect to $\Gamma_2$. 
Thus, if the Euler number of $M_4$ is divisible by $24$,
then the automorphic form $E(M_4;\,Z)$ has trivial character.
The Fourier expansion of the cusp forms  
$\Delta_2(Z)$, $\Delta_{11}(Z)$ and 
$\frac {\Delta_{11}(Z)}{\Delta_{2}(Z)}$
coincide with  the Weyl--Kac--Borcherds denominator formula
of  generalized Kac--Moody super-algebras 
with generalized Cartan matrix $A_{1,II}$, 
$A_{2,II}$ and $A_{2,0}$ respectively:
$$
A_{2,II}=
\left(\smallmatrix
\hphantom{-}{2}&{-2}&{-6}&{-2}\cr
{-2}&\hphantom{-}{2}&{-2}&{-6}\cr
{-6}&{-2}&\hphantom{-}{2}&{-2}\cr
{-2}&{-6}&{-2}&\hphantom{-}{2}\cr
\endsmallmatrix\right),
\quad
A_{2,0}=
\left(\smallmatrix\hphantom{-}2&-2&-2\\
-2&\hphantom{-}2&\hphantom{-}0\\
-2&\hphantom{-}0&\hphantom{-}2
\endsmallmatrix\right),
\quad
A_{2,I}=
\left(\smallmatrix
\hphantom{-}{2}&{-2}&{-4}&\hphantom{-}{0}\cr
{-2}&\hphantom{-}{2}&\hphantom{-}{0}&{-4}\cr
{-4}&\hphantom{-}{0}&\hphantom{-}{2}&{-2}\cr
\hphantom{-}{0}&{-4}&{-2}&\hphantom{-}{2}\cr
\endsmallmatrix\right)
$$
(see \cite{GN1}--\cite{GN4}).
Thus, the formula \thetag{3.4} gives us three particular cases
of Calabi--Yau $4$-folds of Kac--Moody type when 
the second quantized elliptic genus is a power of 
the denominator function of the corresponding  Lorentzian Kac--Moody
algebra:
$$
\align
E(M_4;\,Z)&=\Delta_{11}(Z)^{-\chi_0}
\qquad\qquad\text{if } \chi_1=0\\ 
E(M_4;\,Z)&=\biggl(\frac{\Delta_{11}(Z)}{\Delta_{2}(Z)}\biggr)^{-\chi_0} 
\qquad\text{if } \chi_0(M)=-\chi_1(M)\\
E(M_4;\,Z)&=\Delta_{2}(Z)^{\chi_1}\qquad\qquad\quad\text{if } \chi_0(M)=0.
\endalign 
$$

At this point  we would like to discuss some 
relations of Hecke type between the basic Siegel modular
forms occurred in SQEG-construction. It would be interesting 
to find a geometric or  a physical interpretation of such relations.
We recall that  the Hecke operators
$T_-(m)$ on the space of Jacobi forms of weight zero 
play the main role in the construction of the Borcherds products
\thetag{3.3} (see \cite{B2}, \cite{GN1}, \cite{GN3})
and in the proof of the identity  \thetag{3.1}
for the second quantized elliptic genus (see \cite{DMVV}).
We start with a relation
$$
\psi_{0,2}^{(2)}=\phi_{0,1}|T_-(2)-2\phi_{0,2}.
$$
Thus, according to \cite{GN4, Theorem 3.3} we can represent 
$\Delta_{11}(Z)$ in the following form
$$
\Delta_{11}(Z)=
\Delta_5(Z)
\Delta_5(
\left(\smallmatrix \tau&2z\\ 2z&4\omega \endsmallmatrix\right)
)
\Delta_5(
\left(\smallmatrix\tau&z\\z&{\omega+\frac 1{2}}\endsmallmatrix\right)
)
\Delta_2(Z)^{-2}.
$$
Moreover we can write  $\Delta_{11}^2$ only in terms of 
the cusp form $\Delta_2$  since
$$
2\psi_{0,2}^{(2)}=\phi_{0,2}|T_0(2),
$$
where $T_0(2)$ is another Hecke operator which does not change
the index of Jacobi forms. If $f(n,l)=g(8n-l^2)$ is the Fourier
coefficient of  $\phi_{0,2}$ (for a prime index $t$
Fourier coefficients depend only on the norm $4tn-l^2$), 
then the Fourier coefficient
 $f_2(n,l)=g_2(8n-l^2)$ of  $\phi_{0,2}|T_0(2)$ is given by 
the formula
$$
g_2(N)=8g(4N)+ 2\biggl(\frac {-N}2\biggr)g(N)+
g\biggl(\frac {N}{4}\biggr).
$$
According \cite{GN2, Theorem A.7},
$\Delta_{11}(Z)^2={[\Delta_2(Z)]_{T(2)}}/{\Delta_2(Z)^4}$,
where 
$$
\align
[\Delta_2(Z)]_{T(2)}=
\kern-10pt\prod_{a,b,c\,mod\, 2}&\kern-8pt
\Delta_2 ({\tsize \frac{z_1+a}2 ,
\tsize\frac{z_2+b}2, \tsize\frac{z_3+c}2})
\prod_{a \,mod\, 2}
\Delta_2 ({\tsize\frac{z_1+a}2,z_2,2z_3})\,
\Delta_2 ({\tsize 2z_1, z_2, \frac{z_3+a}2})\\
\vspace{2\jot}
{}&\times\Delta_2 ({\tsize 2z_1, 2z_2, 2z_3})
\prod_{b\,mod\, 2}
\Delta_2 ({\tsize 2z_1, -z_1+z_2, \frac{z_1-2z_2+z_3+b}2}).
\endalign
$$
\smallskip
{\bf 3. The case of $\bold C\bold Y_{\bold 6}$.}
In this case the elliptic genus of a Calabi--Yau $6$-fold
is a sum of three basic Jacobi forms
$\psi_{0,3}^{(i)}$ (see \thetag{1.11}).
Thus  there are three  basic Siegel modular forms
$$
\align
\Delta_1(Z)&=\ml(\psi_{0,3}^{(1)})=
\hbox{Lift}(\eta(\tau)\vartheta(\tau,z))\\
{}&=\sum_{M\ge 1}
\sum
\Sb
m >0,\,l\in \bz\\
\vspace{0.5\jot} n,\,m\equiv 1\,mod\,6\\
\vspace{0.5\jot}
4nm-3l^2=M^2
\endSb
\hskip-4pt
\biggl(\dsize\frac{-4}{l}\biggr)
\biggl(\dsize\frac{12}{M}\biggr)
\sum\Sb a|(n,l,m)\endSb \biggl(\dsize\frac{6}{a}\biggr)
q^{n/6}y^{l/2}s^{m/2},\\
D_6(Z)&=\ml(\psi_{0,3}^{(2)})=
\hbox{Lift}\biggl(\eta(\tau)^{12}\phi_{0,\frac{3}2}
\biggr), \quad
\Delta_{17}(Z)=\ml(\psi_{0,3}^{(3)})
\endalign
$$
with divisors
$H_1$, $H_4$ and $H_1+H_9$ respectively.
(These Humbert modular surfaces  have only one component in $\Cal A_6^+$).
The Fourier expansion of the cusp forms $\Delta_1$, $D_6$ and 
$\Delta_7(Z)=\Delta_{1}(Z)D_6(Z)$ determine
three  generalized Kac--Moody Lie super-algebras
with the following  Cartan matrix of the real simple roots: 
$$
A_{3,II}=
\left(\smallmatrix
\hphantom{-}{2}&{-2}&{-10}&{-14}&{-10}&{-2}\cr
{-2}&\hphantom{-}{2}&{-2}&{-10}&{-14}&{-10}\cr
{-10}&{-2}&\hphantom{-}{2}&{-2}&{-10}&{-14}\cr
{-14}&{-10}&{-2}&\hphantom{-}{2}&{-2}&{-10}\cr
{-10}&{-14}&{-10}&{-2}&\hphantom{-}{2}&{-2}\cr
{-2}&{-10}&{-14}&{-10}&{-2}&\hphantom{-}{2}\cr
\endsmallmatrix\right),\ 
A_{3,I}=
\left(\smallmatrix
\hphantom{-}{2}&{-2}&{-5}&{-1}\cr
{-2}&\hphantom{-}{2}&{-1}&{-5}\cr
{-5}&{-1}&\hphantom{-}{2}&{-2}\cr
{-1}&{-5}&{-2}&\hphantom{-}{2}\cr
\endsmallmatrix\right),\ 
A_{3,0}=
\left(\smallmatrix
\hphantom{-}{2}&{-2}&{-2}\cr
{-2}&\hphantom{-}{2}&{-1}\cr
{-2}&{-1}&\hphantom{-}{2}\cr
\endsmallmatrix\right).
$$
It follows that  SQEG of an arbitrary Calabi--Yau $6$-fold 
can be represented as

$$
E(M_6, Z)=\Delta_{17}(Z)^{-\chi_0}\Delta_7(Z)^{\chi_1}\Delta_1(Z)^{-\chi_2}.
$$
{\it Thus there are three types of  manifolds of Kac--Moody type}:
($\chi_{0}=0$, $\chi_1=0$) or  
$\chi_2=-\chi_1$ or $\chi_2=0$.
Similar to $\Delta_{11}$ the cusp form $\Delta_{17}$ can be represented 
using  a Hecke type  product of $\Delta_5$
(see \cite{GN4, Theorem 3.3}):
$$
\Delta_{17}(Z)=
\biggl[\Delta_5(\left(\smallmatrix\tau&3z\\3z&9\omega\endsmallmatrix\right))
\prod_{b\,mod\,3} 
\Delta_5(\left(\smallmatrix\tau&z\\z&\omega+\frac b{3}\endsmallmatrix\right))
\biggr] \Delta_1(Z)^{-3}.
$$

\smallskip
{\bf 4. The case of $\bold C\bold Y_{\bold 8}$ and 
$\bold C\bold Y_{\bold 10}$.}
In the case  $d=8$ we have four basic weak Jacobi forms of index $4$
(see \thetag{1.9}).
The exponential lifting of $\phi_{0,4}(\tau,z)$ 
is the ``most odd" even Siegel 
theta-constant 
$$
\ml(\phi_{0,4})(Z)=\Delta_{1/2}(Z)=
\frac{1}2\sum\Sb n,\,m\in \bz \endSb
\,\biggl(\dsize\frac{-4}{n}\biggr)\biggl(\dsize\frac{-4}{m}\biggr)
q^{n^2/8}y^{nm/4}s^{m^2/8}
$$
which is the  denominator function of a generalized 
Lorentzian Kac--Moody super-algebra with a real root system of
{\it  parabolic 
type} (see \cite{GN2--GN3}). 
For arbitrary Calabi--Yau $8$-fold we have (see \thetag{1.11}):
$$
E(CY_8;\, Z)=\Delta_{1/2}(Z)^{\chi_1-\chi_3}\cdot
\Delta_{5}(\left(\smallmatrix \tau&2z\\2z&4\omega\endsmallmatrix\right))
^{-\chi_2}\cdot \Delta_{12}(Z)^{\chi_1}\cdot
\Delta_{23}(Z)^{-\chi_0},
$$
where
$$
\Delta_{12}(Z)=\ml(\psi^{(3)}_{0,4})\in \frak M_{12}(\Gamma_4^+),\quad
\Delta_{23}(Z)=\ml(\psi^{(4)}_{0,4})\in \frak M_{23}(\Gamma_4^+)
$$
are $\Gamma_4^+$-modular forms with trivial character.
We remark that  all forms  above have multiplicities one along divisors!
From the formulae above it follows that for
$d=2$, $4$, $6$, $8$ the divisor of $E(M_d; Z)$
is determined by the divisor of Hodge anomaly of $M_d$.
For $d=10$ this is not the case. From \thetag{1.12} we see
that for $d=10$ one new divisor appears.
This is the Humbert surface
$H_5=\pi_5^+\{\tau+5z+5\omega=0\}\subset \Cal A_5^+$.
The second quantized elliptic genus $E(M_{11}; Z)$
{\it  is anti-invariant with respect to the involution
defined by this rational quadratic divisor}.

We remark also that for $d<12$ (or $d=13$) we have
$$
\align
E(CY_d;\, Z)&=E(CY_d^{mir};\, Z)^{-1}\qquad\  (\text{if } d \text{ is odd})\\
E(CY_d;\, Z)&=E(CY_d^{mir};\, Z)\qquad\quad
 \ \, (\text{if } d \text{ is even})
\endalign
$$
where  $CY_d^{mir}$ is the mirror partner of Calabi--Yau $d$-fold.
\subhead
3.3. SQEG of Calabi--Yau $3$-folds and $5$-folds
\endsubhead
The elliptic genus of  a Calabi--Yau $3$-folds $M_3$ 
is defined uniquely up to a constant (see Example 1.5):
$$
\chi(M_3,\tau,z)=\frac 1{2}{e(M_3)}\phi_{0, \frac 3{2}}(\tau,z).
$$ 
According \thetag{3.3} and Theorem 3.1 the product
$$
\multline
\Phi_3 (Z)=\ml(\phi_{0,\frac{3}2}(\tau,2z))=
\\
q^{\frac{1}{12}}y^{\frac{1}2}s^{\frac{1}2}
\prod\Sb n,l,m\in \bz\\
\vspace{0.5\jot}
(n,l,m)>0\endSb
 (1-q^ny^ls^{tm})^{c(nm,l)}
\in \frak M_0(\Gamma_6^+, v_{12})
\endmultline
\tag{3.5}
$$
is an automorphic function of weight $0$ with respect to 
$\Gamma_6^+$ with a character of order $12$. (This character is induced
by the square $v_\eta^2$ of the multiplier system of the Dedekind eta-function
$\eta(\tau)$ and by the binary character $v_H$ of the  Heisenberg group
defined in \S 1.) Moreover 
$$
\hbox{div}_{\Cal A_6^+}(\Phi_3 (Z))=H_1(0)-H_1(5)=
\pi_6^+(\{Z\in \bh_2\,|\,z=0\})-
\pi_6^+(\{Z\in \bh_2\,|\,\tau+5z+6\omega=0\}).
$$
The divisor $H_1(5)$ appears because
the Fourier expansion of $\phi_{0, \frac 3{2}}(\tau,2z)$ 
contains the term  $-qy^5$.
We remark that $\Phi_3 (Z)$ is anti-invariant with respect to the reflection
in the vectors which define the Humbert surfaces 
$H_1(0)$ and $H_1(5)$.
The character of the automorphic function $\Phi_3 (Z)$ has 
the maximal possible  order $12$ since
$\Gamma_{6}^+/(\Gamma_{6}^+)^{com}=\bz/2\times \bz/12$
(see \cite{GH2, Theorem 2.1}). 
One can say that $\Phi_3 (Z)$ {\it is a Siegel modular function
with the simplest possible divisor}.
Such a function should have a zero and pole inside $\bh_2$.
The simplest divisor is a rational quadratic divisor
(a Humbert surface) with discriminant  one.
The case $t=6$ is the first case when there are two 
divisors of discriminant one. 
$\Phi_3 (Z)$ has zero and poles of order one along these 
surfaces. Thus {\it we can consider
 $\Phi_3 (Z)$ as the best analogue of the modular invariant 
$j(\tau)$ in the case of Siegel modular functions!}

It is interesting that  $\Phi_3 (Z)$
is related to Siegel modular threefolds of geometric genus 
$1$ and $2$. More exactly, the square of  $\Phi_3 (Z)$
can be written as quotient of two $\Gamma_6^+$-cusp forms
constructed in  \cite{GN2, Example 4.6}
$$
\Phi_3 (Z)^2=
\frac{\hbox{Lift}\bigl(
\eta(\tau)^3\vartheta(\tau,z)^2\vartheta(\tau,2z)\bigr)}
{\hbox{Lift}(\eta(\tau)^5\vartheta(\tau,2z))}=
\frac{\hbox{Exp-Lift}(3\phi_{0,3}^2-2\phi_{0,2}\phi_{0,4})}
{\hbox{Exp-Lift}(5\phi_{0,3}^2-4\phi_{0,2}\phi_{0,4})}.
$$
The modular form in the numerator 
$
F_1(Z)=\hbox{Lift}\bigl(
\eta(\tau)^3\vartheta(\tau,z)^2\vartheta(\tau,2z)\bigr)
\in \frak N_3(\Gamma_6, v_2)
$
is  cusp form of weight $3$ with divisor $3H_1(1)+2H_1(5)+H_4$.
Let us introduce the subgroup
$\Gamma_6(v_2)=\hbox{Ker}(v_2)\subset \Gamma_6$ of index $2$
($v_2$ is character of order $2$) and the covering
$$
\Cal A_6(v_2)=\Gamma_6(v_2)\setminus \bh_2
\overset{2:1}\to \longrightarrow {\Cal A}_6.
$$
It is known that $h^{3,0}(\Cal A_6(v_2))=1$ and the cusp form
$$
F_1(Z)d\tau\wedge dz\wedge d\omega\in H^{3,0}(\Cal A_6^{o}(v_2), \bc)
$$ 
defines the unique canonical differential form on a smooth
model of $\Cal A_6(v_2)$ (see \cite{GH2, Theorem 3.1}).
The form in the numerator
$F_2(Z)={\hbox{Lift}(\eta(\tau)^5\vartheta(\tau,2z))}
\in \frak N_3(\Gamma_6, v_3)$
is a cusp form of weight $3$ with a character of order $3$.
The Siegel threefold
$$
\Cal A_6(v_3)=\Gamma_6(v_3)\setminus \bh_2
\overset{3:1}\to \longrightarrow {\Cal A}_6
$$
has $h^{3,0}=2$ and $F_2(Z)$ defines one of  
the canonical differential forms on it.

Let us consider an arbitrary  Calabi--Yau $3$-fold $M_3$ 
with Hodge numbers  $h^{1,1}$ and $h^{2,1}$.
Its modified  SQEG is the modular function
$$
E^{(2)}(M_3;\,Z)=\Phi_3 (Z)^{h^{2,1}-h^{1,1}}.
$$
In particular $\Phi_3 (Z)$ is the modified SQEG of a Calabi--Yau
$3$-fold with Euler number $-2$.
The first divisor 
$(h^{2,1}-h^{2,2})H_1(0)$ of $E^{(2)}(M_3;\,Z)$
comes from the Hodge anomaly $H(M_3;\,Z)$,  
the second $(h^{2,2}-h^{2,1})H_1(5)$ is the additional divisor
which is, in some sense, ``mirror symmetric" to the first one, because if
$M_3^{mir}$ is a mirror partner of $M_3$, then
$$
E^{(2)}(M_3^{mir};\,Z)=E^{(2)}(M_3;\,Z)^{-1}.
$$
It turns out (see \cite{GN5}) that  {\it the functions 
$\Phi_3 (Z)^{\pm 1}$ determine Lorentzian Kac--Moody super-algebras
with an infinite system of real simple roots of the third possible type}.
This is the  so-called {\it hyperbolic type} when the infinite system 
of real simple roots  has a  ``limit" line
(see \cite{N}).
The Lorentzian Kac--Moody algebras related 
to SQEG$(CY_d)$ ($d=2,\,4,\,6,\,8$)  have elliptic type
(the system of real  simple roots  is finite) or parabolic type
(the system of real simple roots has geometry similar
to the real simple roots of the fake monster Lie algebra constructed
by Borcherds in \cite{B1}). 

\smallskip
In the case  $d=5$ the elliptic genus is  again defined uniquely
by the Euler number.
The basic Jacobi form is 
$
\phi_{0, \frac{5}2}=\phi_{0, \frac{3}2}\phi_{0,1}
$
(see Example 1.5). We have that
$$
\Phi_5 (Z)=\ml(\phi_{0,\frac{5}2}(\tau, 2z))
\in \frak M_0(\Gamma_{10}^+)
$$
is  an automorphic function of weight $0$ with respect to 
$\Gamma_{10}^+$ with trivial character.
The automorphic function $\Phi_5 (Z)$ has the divisor
consisting of four irreducible components
$$
\hbox{div}_{\Cal A_{10}^+}\bigl(\Phi(Z)\bigr)=
H_9(3)-H_9(7)+12H_{1}(1)-12H_{1}(9).
$$
For arbitrary Calabi--Yau fivefold $M_5$ with Euler number $\chi(M_5)$
we have
$$
E^{(2)}(M_5;\,Z)=\Phi_5(Z)^{-e(M_5)/{24}}.
$$


\Refs 
\widestnumber\key{DMVV}

\ref
\key AYS
\by O. Aharony, S. Yankielowicz, A.N. Schellekens
\paper 
Charge sum rules in $N=2$ theories
\jour Nucl. Phys. 
\vol B418 
\yr 1994
\pages 157
\endref 



\ref
\key B1
\by R. Borcherds
\paper The monstrous moonshine and monstrous Lie superalgebras
\jour Invent. Math.
\vol 109
\yr 1992
\pages 405--444
\endref


\ref
\key B2
\by R. Borcherds
\paper Automorphic forms on $O_{s+2,2}$ and
infinite products
\jour Invent. Math. \vol 120
\yr 1995
\pages 161--213
\endref

\ref
\key D 
\by R. Dijkgraaf
\paper The Mathematics of Fivebranes
\jour Documenta Mathem. ICM-1998
\yr 1998 
\endref 


\ref
\key DVV 
\by R. Dijkgraaf, E. Verlinde and H. Verlinde 
\paper Counting dyons in $N=4$ string theory 
\jour Nucl. Phys.
\vol  B484 (1997)
\yr 1997
\pages 543--561  
\endref 

\ref
\key DMVV 
\by R. Dijkgraaf, G. Moore,  E. Verlinde, H. Verlinde 
\paper Elliptic genera of symmetric products
and second quantized strings
\jour Commun. Math. Phys.
\vol 185
\yr 1997
\pages 197--209
\endref 

\ref
\key
EOTY
\by
 T. Eguchi, H. Ooguri, A. Taormina, S.-K. Yang,
\paper
Superconformal Algebras and String Compactification on Manifolds
with $SU(N)$ Holonomy
\jour
 Nucl. Phys. 
\yr 1989
\vol  B315
\pages 193--221
\endref


\ref
\key G\"o
\by
L. G\"ottsche
\paper The Betti numbers of the Hilbert Scheme of Points
on a Smooth Projective Surface
\jour Math. Ann. 
\vol 286
\yr 1990
\pages 193--297
\endref



\ref
\key G1
\by V\. Gritsenko
\paper Modular forms and moduli spaces of Abelian and K3 surfaces
\jour Algebra i Analyz
\vol 6:6
\yr 1994
\pages 65--102
\transl\nofrills  English transl. in
\jour St.Petersburg Math. Jour.
\vol 6:6
\yr 1995
\pages 1179--1208
\endref

\ref
\key G2
\by V\. Gritsenko
\paper Modified Witten genus
\jour Preprint MPI
\yr 1999
\endref

\ref
\key GH1
\by V. Gritsenko, K. Hulek
\paper Minimal Siegel modular threefolds
\jour Mathem. Proc. Cambridge Phil. Soc. 
\vol 123
\yr 1998
\pages  461--485
\endref



\ref
\key GH2
\by V. Gritsenko, K. Hulek
\paper Commutator coverings of Siegel threefolds
\jour  Duke Math. J.
\vol 94
\yr 1998
\pages 509--542
\endref




\ref
\key GN1
\by V.A. Gritsenko, V.V. Nikulin
\paper Siegel automorphic form correction of some Lorentzi\-an
Kac--Moody Lie algebras
\jour Amer. J. Math.
\yr 1997 
\vol 119
\pages  181--224
\endref



\ref
\key GN2 
\by V.A. Gritsenko, V.V. Nikulin
\paper The Igusa modular forms and ``the simplest''
Lorentzian Kac--Moody algebras
\jour Matem. Sbornik 
\yr 1996 \vol 187 
\pages 1601--1643 
\endref

\ref
\key GN3
\by V.A. Gritsenko, V.V. Nikulin
\paper Automorphic forms and Lorentzian Kac-Moody algebras.
Part I 
\jour International  J. of Mathem.
\vol 9
\yr 1998
\pages 153--199
\endref

\ref
\key GN4 
\by V.A. Gritsenko, V.V. Nikulin
\paper Automorphic forms and Lorentzian Kac-Moody algebras.
Part II 
\jour International  J. of Mathem.
\vol 9
\yr 1998
\pages 201--275
\endref

\ref
\key GN5
\by V.A. Gritsenko, V.V. Nikulin
\paper Second quantized elliptic genus 
of Calabi--Yau  $3$-folds 
\jour Preprint (in preparation) 
\endref


\ref
\key HM1 
\by J. Harvey, G. Moore 
\paper Algebras, BPS-states, and strings 
\jour Nucl. Physics. 
\vol B463 
\yr 1996
\pages 315--368  
\endref 


\ref
\key H1
\by F. Hirzebruch
\paper Elliptic genera of level N for complex manifolds
\inbook Differential geometrical Methods
in Theoretical Physics 
\eds K. Bleuler, M. Werner
\publ Kluwer Acad. Publ.
\yr 1988
\pages 37--63
\moreref Appendix III to [HBJ]
\endref

\ref
\key H2
\by F. Hirzebruch
\paper Letter to V. Gritsenko from  11 August 1997
\endref


\ref
\key HBJ
\by F. Hirzebruch, T. Berger, R. Jung
\book Manifolds and Modular forms
\publ Aspects of Math. {\bf E20}
Vieweg--Verlag
\yr 1992
\endref


\ref
\key HH
\by F. Hirzebruch, T. H\"ofer
\paper
On the Euler Number of an Orbifold
\jour
Math. Ann. 
\yr 1990
\vol 286
\pages 255--260
\endref


\ref
\key H\"o
\by G. H\"ohn
\paper  Komplex elliptische Geschlechter und 
$S^1$-\"aquivariante Kobordismustheorie
\jour
 Diplomarbeit
\yr 1991
\publaddr Bonn
\endref


\ref
\key K
\by I. Krichever
\paper Generalized elliptic genera and Baker--Akhiezer functions
\jour Mat. Zametki
\vol 47
\yr 1990
\pages 34--45
\endref

\ref
\key KYY
\by
T. Kawai, Y. Yamada, S.-K. Yang
\paper
Elliptic Genera
and N=2 Superconformal Field Theory
\jour
Nucl. Phys. 
\vol B414
\yr 1994
\pages 191-212
\endref


\ref
\key L
\by P.S. Landweber Ed.
\book Elliptic Curves and Modular Forms in 
Algebraic Topology
\publ Springer-Verlag
\yr 1988
\endref 

\ref
\key LW
\by A. S. Libgober, J.W. Wood
\paper Uniqueness of the complex structure of K\"ahler manifolds
of certain homotopy types
\jour J. Differential Geometry
\vol 32 \yr 1990
\pages 139--154
\endref

\ref
\key M
\by
G. Moore
\paper
String duality, automorphic forms and generalized 
Kac--Moody algberas
\jour
Nucl. Phys. Proc. Suppl.
\vol 67
\yr 1998
\pages 56--67
\endref


\ref
\key N
\by V. V. Nikulin
\paper K3 surfaces with interesting groups of automorphisms
\jour  to appear in J. Math. Sciences
\moreref alg-geom/9701011 
\endref


\ref
\key O
\by S. Ochanine 
\paper Signature modulo 16, invariants de Kervaire g\'en\'eralis\'es,
et nombres caract\'eris\-tiques dans la K-th\'eorie r\'eelle
\jour M\'emoire de la Soc. Math. de France; Nouvelle S\'eries
\vol 5 \yr 1981
\endref

\ref
\key R
\by V. A. Rokhlin
\paper New results in the theory of the varieties of dimension 4
\jour Doklady AN USSR
\vol 84 \yr 1952
\pages 221--224
\endref


\ref
\key S
\by S. M. Salamon
\paper On the cohomology of K\"ahler and hyper-K\"ahler manifolds
\jour Topology
\vol 35 \yr 1996
\pages 137--155
\endref

\ref
\key SVW
\by
S. Sethi, C. Vafa, E. Witten
\paper Constraints on low-dimensional string
compactifications
\jour   Nucl. Phys. 
\vol 
B480
\yr 1996
\pages 213--224
\endref

\ref
\key W1
\by
E. Witten
\paper Elliptic genera and Quantum Field Theory
\jour Commun. Math. Phys. 
\vol 109
\yr 1987
\pages 525--536
\endref

\ref
\key W2
\by E. Witten
\paper The index of the Dirac operator in loop space
\pages 161--181
\moreref in [L]
\endref

\ref
\key Z
\by
D. Zagier
\paper Note on the Landweber--Strong elliptic genus 
\pages 216--224
moreref in [L]
\endref

\endRefs
\enddocument
\end